\documentclass{article}

\usepackage{amsmath,amsthm,amsfonts,amssymb}

\theoremstyle{plain}
\newtheorem{theorem}{Theorem}[section]
\newtheorem{lemma}[theorem]{Lemma}

\newtheorem{corollary}[theorem]{Corollary}
\newtheorem{proposition}[theorem]{Proposition}

\theoremstyle{definition}
\newtheorem{definition}[theorem]{Definition}

\theoremstyle{remark}

\newtheorem{remark}[theorem]{Remark}
\newtheorem{example}[theorem]{Example}

\DeclareMathOperator*{\colim}{colim}

\usepackage{mathrsfs}
\usepackage{amsfonts}
\usepackage[all]{xy}
\usepackage{fancyhdr}
\usepackage[left=1in,right=1in,top=1in,bottom=1in]{geometry}

\def \C              {\mathcal C}
\def \P              {\mathbb{P}}
\def \K              {\mathbb{K}}
\def \F              {\mathbb{F}}
\def \R              {\mathbb{R}}
\def \Lf             {\mathbb{L}}
\def \N              {{\mathcal N}_{ {\mathbb{K}}}}
\def \L              {{\mathcal L}_{ {\mathbb{K}}}}
\def \A              {{\mathcal Ass}_{ {\mathbb{K}}}}
\def \CH             {{\mathcal CH}_{ {\mathbb{K}}}}
\def \T              {\mathbb{T}_1}

\def \g              {\bf g}
\def \h              {\bf h}

\begin{document}

\title{Functional Equations and Their Related Operads}
\author{Vahagn Minasian}
\date{}

\maketitle

\begin{abstract}
Using functional equations, we define functors that generalize standard examples 
from calculus of one variable. Examples of such functors are discussed and their Taylor 
towers are computed. We also show that these functors factor through objects enriched 
over the homology of little $n$-cubes operads and discuss the relationship between 
functors defined via functional equations and operads. In addition, we compute the 
differentials of the forgetful functor from the category of $n$-Poisson algebras in 
terms of the homology of configuration spaces.  

\vspace{10pt}

\noindent
MCS: 55U15 (18D50, 55P99)
\end{abstract}

\section*{Introduction}

In late eighties, Tom Goodwillie in a series of papers 
(~\cite{Good1},~\cite{Good2},~\cite{Good3}) introduced a construction that plays the 
role of traditional Taylor series for functors of spaces. More precisely, for the given 
functor $F$, a tower of functors $\{P_n F\}$ was constructed with $P_n F$ of degree 
$n$. Moreover, this tower (when evaluated at a space $X$ satisfying certain connectivity 
conditions) approximates $F$ in some appropriate sense. Since then 
Goodwillie Calculus has been used extensively and new methods were developed 
to provide alternative descriptions in different categories. 
In particular, in~\cite{Randy}, B.Johnson and R.McCarthy introduce a construction using 
cotriples that produces a Taylor tower of functors from any pointed category with 
coproducts to the category of chain complexes $Ch(\K)$ over a commutative ring $\K$. 
     
Traditionally, after introducing the Taylor series for a function of one variable, one 
studies the Taylor series of simple functions, such as the exponential and logarithmic 
functions or the function $\frac{x}{1-x}$, since the derivatives of these functions can be 
easily computed. Not unexpectedly, computing Taylor towers of functors is more complicated. 
However some of the techniques from traditional real variable calculus apply. For
example, if a one variable function satisfies a sufficiently simple functional equation,
its Taylor series can often be computed.  

In~\cite{Randy}, the authors observed that exponential functors are defined as a solution
to a rather elementary functional equation, which allowed them to compute the cross effects
and consequently, the layers (or fibers) of Taylor towers of these functors.  

We follow this approach and define new classes of functors via functional equations.
In addition to exponential functors, these are functors of type $\frac{x}{1-x}$,
logarithmic functors and functors of type $f_n$, which mimic the one variable 
functions $\frac{x}{1-x}$, $-log(1-x)$, and 
$f_n=(1- \frac{x}{2^{n-1}})^{-2^{n-1}}-1$ respectively. 

A crucial observation is that forgetful functors from the category of 
(non-unital) commutative algebras, (non-unital) associative algebras, Lie 
algebras and $n$-Poisson algebras to chain complexes over $\K$ are examples 
of exponential functors, functors of type $\frac{x}{1-x}$,
logarithmic functors and functors of type $f_n$ respectively.
This leads to some interesting applications of our computations. For example
in~\cite{Minas}, we 
follow the approach suggested in~\cite{Randy} and 
use the Taylor towers of exponential functors to 
derive the topological analogue (in framework of~\cite{EKMM})
of the commutative to associative spectral sequence,
also known as the {\it fundamental spectral sequence for} $THH$. 

In a joint work with R.McCarthy (~\cite{McMin}), functors of type $\frac{x}{1-x}$
are employed to develop a set up (analogous to that of Quillen in~\cite{Quillen} for 
discrete case), which allows a discussion on homology and abelianization functors 
for $A_\infty$ algebras. 
Furthermore, specializing the results of (~\cite{McMin}) to functors of type 
$\frac{x}{1-x}$ or type $f_n$ and their towers, gives us some splitting
results for associative and $n$-Poisson algebras.       

On a more curious point, this approach of associating algebraic structures 
with real variable functions gives an interesting angle of looking at some 
classical algebraic results. For example, the simple equation
$$e ^{-\log(1-x)} -1 = \frac{x}{1-x}$$
is at the root of the Poincar\'e-Birkhoff-Witt Theorem. 

In other words, studying these functors of various types aids in understanding 
different algebraic structures, which suggests a relationship between such functors
and some of the standard operads. More precisely, we construct a sequence of 
functors (defined as solutions to certain functional equations), namely functors 
of type $\frac{x}{1-x}$, $f_2, f_3, \cdots$, and $exponential$ functors, which have 
the same set of algebras of different types associated with them as the 
homology of little $n$-cubes operads $\{ e_n \}_{n=1} ^{\infty}$, i.e. the 
associative algebras (for $e_1$ and functors of type $\frac{x}{1-x}$), 
$n$-Poisson algebras (for $e_n$ and type $f_n$ with $n \geq 2$) and commutative 
algebras (for $e_{\infty}$ and $exponential$ functors). 

In addition, to understanding the functors of different types defined via functional 
equations, the exploration of this relationship with operads is the primary goal of 
this work.      

The paper is structured as follows. 

In Sections 1, 2, and 3 we set up the groundwork, as well as briefly recall some of the 
basic constructions from Goodwillie Calculus, little $n$-cubes operads and the theory 
of cotriples, triples and algebras associated with them. References for in depth study 
of this objects are provided.  

The rest of the paper is built with analytic functions of one variable in mind. In other 
words, the real variable function which is to be imitated, is introduced, then the 
corresponding functor is defined, followed by a discussion of the algebraic structures it 
carries and examples. Thus, Section 4 focuses on exponential functors, Section 5 introduces
functors of type $\frac{x}{1-x}$, in Section 6 the Taylor towers of these functors 
are computed. In Section 7, we discuss logarithmic functors, as well as compute their Taylor
towers. Section 8 deals with functors of type $f_n$, whose Taylor towers are computed in 
Section 9. 
Section 10 considers the set of different types 
as a whole. In addition, the relationship with the homology of little $n$-cubes operads is 
explored in this section. In particular, the forgetful functor from the category 
of $n$-Poisson algebras to chain complexes is considered in light of F.Cohen's computations
of the homology of configuration spaces.   

{\bf Acknowledgments.} The author would like to thank Randy McCarthy for his encouragement
and many helpful comments.    
 
\section{Triples, Cotriples and Algebras}

Triples arise in category theory as an abstraction of the notion of an algebraic
structure on an underlying space. Let $\C$ be a pointed category with coproducts.

\begin{definition}
A triple $(T, \mu, \eta)$ on $\C$ is a functor $T:\C \to \C$, together with 
natural transformations $\mu :TT \to T$ and $\eta : Id \to T$, such that 
the following diagrams commute: 

\hspace{60pt}
$
\xymatrix{
TTT \ar[r]^{T \mu} \ar[d]
&
TT \ar[d]^{\mu}\\
TT \ar[r]^{\mu}
&
T}
$
\hspace{60pt}
$
\xymatrix{
T \ar[r]^-{T \eta} \ar[rd]_T
&
TT \ar[d]^\mu 
&
T \ar[l]_-{\eta T} \ar[ld]^T\\
&
T
&}
$

An algebra over a triple is a pair $(X, \rho)$ where $X$ is an object of the category 
$\C$ and $\rho:TX \to X$ is a morphism such that the composition 
$X \stackrel{\eta X}{\to}TX \stackrel{\rho}{\to}X$ is the identity and the obvious 
associativity diagram commutes.  
\end{definition}

Next we recall the definitions of a few common triples. For these examples, let $\C$ be 
a symmetric monoidal category, such as the category of modules over a commutative ring
$\K$ or more generally, the category of chain complexes $Ch(\K)$. The latter category 
is where our subsequent work takes place. Consider the triple
$$\mathbb{T}(X) = \bigotimes X^{\otimes n}.$$  
In the category of modules, the algebras over this triple are precisely the associative 
algebras. Consequently, we will refer to algebras over this triple in any symmetric 
monoidal category as associative algebras. 
Similarly, the triple
$$\P(X) = \bigotimes (X^{\otimes n})/\Sigma_n ,$$
where $\Sigma_n$ is the group of permutations, produces the commutative algebras. 
There are also triples associated to Lie algebras, $n$-Poisson algebras, and other
structured objects.

There is a dual notion of cotriples $(\perp, \epsilon, \delta)$, which consist 
of a functor $\perp:\C \to \C$, together with a natural comultiplication map
$\epsilon:\perp \to \perp \perp$ and a counit $\delta:\perp \to Id$, such that 
the obvious co-associativity and co-unit diagrams commute. 

Adjoint pairs of functors $(F,U)$ provide a good source of cotriples as the composite
$FU$ is always a cotriple. 

\begin{remark}
Let $(\perp, \epsilon, \delta)$ be a cotriple in $\C$. Fix an object $X$ in $\C$. 
Define $\perp^{\ast +1} X$ to be the simplicial object with 
$[n] \mapsto \perp^{n +1} X = \perp \cdots \perp X$, with face and degeneracy maps
$s_i = \perp^{(i)} \delta \perp^{(n-i)}$ and $d_i=\perp^{(i)} \epsilon \perp^{(n-i)}$ .
We will recall this construction when presenting the definition of 
Goodwillie towers. 

\end{remark} 

Chapter 8 of~\cite{Weib} is a good reference for an in depth discussion of triples 
and cotriples. 

\section{Little n-cubes Operads}
\label{sect:operad}

Operads can be defined in any symmetric monoidal category. However, we will only
consider the operads over $Ch(\K)$, the chain complexes over the commutative ring 
$\K$. An operad is a sequence of objects ${\bf a}(k)$, $k \geq 0$, carrying an action
of symmetric groups $\Sigma_k$, with products
$${\bf a}(k) \otimes {\bf a}(j_1)\otimes \cdots \otimes {\bf a}(j_k) \to
  {\bf a}(j_1, \cdots, j_k)$$
which are equivariant and associative in an appropriate sense. For a precise 
definition along with a detailed discussion of operads we refer to ~\cite{Kriz}
or ~\cite{Getz}. 

\begin{example}
Fix a chain complex $V$. Define an operad with 
$$\mathcal{E}_V (k) = Hom(V^{\otimes k}, V).$$
The symmetric group $\Sigma_k$ acts on $\mathcal{E}_V (k)$ through its action on 
$V^{\otimes k}$, and the structure maps are the obvious ones. 
\end{example}

\begin{definition}
An algebra over an operad ${\bf a}$ is a chain complex $A$ together with a morphism
of operads $\rho: {\bf a} \to \mathcal{E}_A$. In other words, $A$ is equipped with 
structure maps 
$$\rho_k:{\bf a}(k)\otimes_{\Sigma_k} A^{\otimes k} \to A$$
which are compatible in a natural sense. 

Equivalently, an ${\bf a}$-algebra can be defined as the algebra over the 
associated triple 
$T_{\bf a} (V) = \bigoplus {\bf a}(k)\otimes_{\Sigma_k} V^{\otimes k}$.
\end{definition}

Next we briefly recall the definition and properties of some of the most important
operads, {\bf the little n-cubes operads}. For details, we again refer to ~\cite{Kriz}
and ~\cite{Getz}. 

Denote by $\F_X (S)$ the configuration space of embeddings of the finite set $S$ 
into the topological space $X$,
$$\F_X (S)=\{x \in X^S| x(s) \neq x(t), s \neq t\}.$$
Abbreviate $\F_{\R ^n}$ by $\F_n$. Define ${\bf e}_n (s)$ to be the homology 
of the configuration space $H (\F_n(S), \K)$, for $s>0$, and ${\bf e}_n (0)=0$.
 ${\bf e}_n (s)$ has an obvious action of $\Sigma_s$, and the sequence 
 $\{{\bf e}_n (s)\}$ assembles into an operad. We call an algebra for the operad
${\bf e}_n (s)$ an $n$-algebra and denote the associated triple by $T_n$. 

Observe that the configuration space $\F_1(k)$ is a disjoint union of contractible 
spaces; the components correspond to permutations $\sigma \in \Sigma_k$ by the rule
$$\sigma \longmapsto \{(x_{\sigma(1)}, \cdots, x_{\sigma(k)}) \in \R^k|
  x_{\sigma(1)}< \cdots < x_{\sigma(k)} \}.$$
It follows that ${\bf e}_n (k)= \K[\Sigma_k]$, and the associated triple is the 
associative algebra triple $T_1 (V) = \bigoplus V^{\otimes k}$.

The configuration spaces of $\R^\infty$ are contractible. Hence the triple associated  
to ${\bf e}_\infty$ is the symmetric algebra triple 
$T_\infty (V) = \bigoplus (V^{\otimes k})_{\Sigma_k}$.


\section{Goodwillie Calculus}
\label{sect:calc}

In this section we give a brief summary of main constructions of Goodwillie Calculus 
as developed by B.Johnson and R.McCarthy. For details we refer to ~\cite{Randy}, which 
is an application of ~\cite{Good1}, ~\cite{Good2} and ~\cite{Good3}. 

We begin by recalling the definition of the cross effects of a functor 
$F: \C \to \mathcal{A}$, where $\C$ is a basepointed category with finite coproducts 
and $\mathcal{A}$ is an abelian category. 
\begin{definition}
The n-th cross effect of $F$ is the functor $cr_n F: \C ^{\times n} \to \mathcal{A}$ 
defined inductively by 

\noindent
$cr_1 F(M) \oplus F( \ast) \cong F(M)$

\noindent
$cr_2 F(M_1, M_2) \oplus cr_1 F(M_1) \oplus cr_1 F(M_2) \cong cr_1 F(M_1 \vee M_2)$

\noindent
and in general,

\noindent
$cr_n F(M_1, \cdots M_n) \oplus cr_{n-1} F(M_1,M_3 \cdots M_n)  
\oplus cr_{n-1} F(M_2,M_3 \cdots M_n)$

\noindent
is equivalent to 
$$cr_{n-1} F(M_1 \vee M_2,M_3 \cdots M_n).$$
\end{definition}

\begin{definition}
Given a functor $F$ from $\C$ to $Ch(\mathcal{A})$, we say that $F$ is degree $n$ is $cr_{n+1} F$
is acyclic as a functor from $\C^{\times n+1}$ to $Ch(\mathcal{A})$. That is, $cr_{n+1} F$ has 
no homology when evaluated on a collection of $n+1$ objects in $\C$. 
\end{definition}

Denote the category of functors of $n+1$ variables from $\C$ to $\mathcal{A}$ that
are reduced in each variable by $Func_\ast (\C^{\times n+1},\mathcal{A})$. Let 
$\bigtriangleup ^\ast$ be the functor from  $Func_\ast (\C^{\times n+1},\mathcal{A})$
to $Func_\ast (\C, \mathcal{A})$ obtained by composing a functor with the diagonal
functor from $\C$ to $\C^{\times n+1}$. The $(n+1)$st cross effect is the right adjoint
to $\bigtriangleup ^\ast$. 

\begin{definition}
Let $\perp _{n+1} =  \bigtriangleup ^\ast \circ cr_{n+1}$ be the cotriple on 
$Func_\ast (\C, \mathcal{A})$ obtained from the adjoint pair 
$(\bigtriangleup ^\ast, cr_{n+1})$. We define $P_n$ to be the functor from 
$Func_\ast (\C, \mathcal{A})$ to $Func_\ast (\C, Ch_{\geq 0} \mathcal{A})$ given by
$$P_n F(X) = hocofiber [N(\perp_{n+1} ^{\ast +1}) \stackrel{\epsilon}{\to} id]$$
where $N$ is the associated normalized chain complex of the simplicial object. 
Furthermore, let $p_n:id \to P_n$ be the natural transformation obtained from the 
mapping cone.
\end{definition}

Next we produce a natural transformation $q_{n}:P_{n} \to P_{n-1}$. 

Observe that we have the following formula relating the $n$'th and $n+1$'st 
cross effects:
\begin{equation}
\label{eq:cross}
cr_{n+1} F(X_1, \cdots, X_{n+1}) = 
cr_2 (cr_n F(X_1, \cdots, X_{n-1},-))(X_n,X_{n+1}).
\end{equation}
Fix an object $X$ of $\C$. Let $G(Y)= cr_n F(X, \cdots, X,Y)$. Then using the fold
map, we have $cr_2 G(X,X) \to G(X \vee X) \to G(X)$, which in turn gives us
$$cr_2 (cr_n F(X, \cdots, X,-))(X,X) \to cr_n F(X, \cdots, X,X).$$
Combine this map with the Equation~\ref{eq:cross} to produce a map
$ cr_{n+1} F(X, \cdots, X) \to cr_{n} F(X, \cdots, X)$, which induces the 
desired map $q_n$. 
By Lemma 2.11 and Theorem 2.12 of ~\cite{Randy} these assemble into a tower
$$
\xymatrix{
&
\cdots \ar[d]^{q_{n+1}}\\
F(X) \ar[r]^{p_n} \ar[rd]^{p_{n-1}}
&
P_n F(X) \ar[d]^{q_{n}}\\
&
P_{n-1} F(X) \ar[d]\\
&
\cdots
}$$ 
such that 

1) The functor $P_n F(X)$ is degree $n$.

2) If $F$ is degree $n$, then $p_n:F \to P_n$ is a quasi-isomorphism. 

3) The pair $(P_n, p_n)$ is universal up to natural quasi-isomorphism with 
respect to degree $n$ functors with natural transformations from $F$.

\begin{definition}
The $n$-th layer or the $n$-th derivative of $F$ is the functor 
$$D_n F(-) (X):= hofiber(q_n)(X).$$ 

\end{definition}

\section{Exponential Functors}

In this section we define the first class of functors that we are interested in 
(exponential functors), and discuss a few of their properties. 
The definition of exponential functors goes back to S.MacLane. However we have learned 
it from~\cite{Randy}, where the authors also prove the following structure 
lemma. Here we present a (different) proof, which is better aligned with the rest 
of this work.  

Reduced functions of real variables (i.e. functions $f(x)$ with  $f(0)=0$) will serve 
as a motivation for our definitions. Observe that functions of the form $f(x)=exp(ax)-1$
satisfy the functional equation 
$$ f(x) f(y) +f(x)+f(y)= f(x+y).$$
Moreover, every reduced  analytic function that satisfies the above equation is of 
the form 
$$f(x)=exp(ax)-1=(exp(x)-1) \circ (ax).$$ 
Thus, we give the following definition. 

\begin{definition}
\label{def:exp}

A reduced functor $F$ from $\C$ to $Ch(\K)$ is exponential if there is a natural 
isomorphism
$$ \alpha _{X,Y}: F(X) \otimes F(Y) \oplus F(X) \oplus F(Y) \rightarrow F(X \vee Y)$$

such that for all objects $X$ and $Y$ of $\C$ the following diagram naturally commutes:
$$
  \xymatrix{
F(X) \otimes F(Y) \oplus F(X) \oplus F(Y) \ar[r]^-{\alpha _{X,Y} }\ar[d]^\cong
&
F(X \vee Y) \ar[d]^\cong \\
F(Y) \otimes F(X) \oplus F(Y) \oplus F(X) \ar[r]^-{\alpha _{Y,X} } &
F(Y \vee X) 
}
$$
\end{definition}

\begin{example}
Let $U_e$ be the forgetful functor from the category $\N$ of non-unital commutative 
differential graded  $\K$-algebras  to the category $Ch(\K)$ of chain complexes. 
Observe that the coproduct of objects $X$ and $Y$ in $\N$ is given by
$$X \otimes Y \oplus X \oplus Y.$$
Thus the functor $U_e$ is exponential. 

\end{example}

In a certain sense $U_e$ is the only example of exponential functors. More precisely, 
continuing with the analogy with functions of real variables, we prove the following 
lemma.

\begin{lemma}

Every exponential functor $F$ from $\C$ to $Ch(\K)$ can be factored as 
$F=U_e \circ \tilde{F}$, where $\tilde{F}$ is a reduced coproduct preserving 
functor from $\C$ to $\N$. 

\end{lemma}

Note that this lemma is the analogue of the earlier statement that every 
reduced analytic function satisfying the exponential equation is of the form 
$(exp(x)-1) \circ (ax)$, as the role of $exp(x)-1$ is played by $U_e$, and the 
reduced coproduct preserving functor $\tilde{F}$ is the appropriate analogue 
of the term $ax$. 

\begin{proof}

Define a product $\mu : F(X) \otimes F(X) \rightarrow F(X)$ by the composite

$$
 \xymatrix{
F(X) \otimes F(X) \ar[r]
& F(X) \otimes F(X) \oplus F(X) \oplus F(X) \ar[r]^-{\alpha _{X,X} }
& F(X \vee X) \ar[r]^-{F(+)}
& F(X)
}
$$

As $\mu$ is commutative (obvious from the definition of an exponential functor), we only 
need to show that that it is associative. 
By naturality, we have a pair of commutative diagrams.

\noindent
$
  \xymatrix{
F(X \vee X) \otimes F(X) \ar[r]^-{F(+) \otimes 1 }\ar[d]^{\alpha_{X \vee X, X}}
&
F(X) \otimes F(X) \ar[d]^{\alpha_{X,X}} \\
F((X \vee X) \vee X) \ar[r]^-{F(+\vee 1)} &
F(X \vee X) 
}
$
\hspace{3pt}
$
  \xymatrix{
F(X) \otimes F(X \vee X) \ar[r]^-{1 \otimes F(+) }\ar[d]^{\alpha_{X, X \vee X}}
&
F(X) \otimes F(X) \ar[d]^{\alpha_{X,X}} \\
F(X \vee (X \vee X)) \ar[r]^-{F(1\vee +)} &
F(X \vee X) 
}
$

\noindent
Thus, it is sufficient to prove that the diagram 
\begin{equation}
\label{expdiag}
  \xymatrix{
F(X) \otimes F(X) \otimes F(X) \ar[rrr] ^-{\alpha_{X,X \vee X} \circ (1 \otimes \alpha_{X,X})}
                               \ar[d]^{\alpha_{X \vee X,X} \circ (\alpha_{X,X} \otimes 1)}
&&&
F(X \vee (X \vee X)) \ar[d]^{F(+)}\\
F((X \vee X) \vee X) \ar[rrr] ^-{F(+)}
&&&
F(X)
}
\end{equation}
commutes.

The following series of equivalences, which follow from the defining equation of exponential 
functors, provides a description of the top horizontal arrow in Diagram~\ref{expdiag}.

$
F(X) \otimes F(X) \otimes F(X) \oplus 3 F(X \vee X) 
\stackrel{1 \oplus 2 \alpha_{X,X} ^{-1} \oplus 1} {\cong} 
$

$
 F(X) \otimes F(X) \otimes F(X) \oplus 2 F(X) \otimes F(X) \oplus 4 F(X) \oplus F(X \vee X)
 \cong
$

$ 
 F(X) \otimes [F(X)\otimes F(X) \oplus F(X) \oplus F(X)] \oplus 4 F(X) \oplus F(X \vee X)
\stackrel{1_{F(X)} \otimes \alpha_{X,X} \oplus 1} {\cong}
$

$ 
 F(X) \otimes F(X \vee X) \oplus 4 F(X) \oplus F(X \vee X) \cong
 F(X \vee X \vee X) \oplus 3 F(X)
$

\noindent
This and an analogous series of equivalences which describes the left vertical map of 
Diagram~\ref{expdiag} are employed to form the following diagram:

$$
  \xymatrix{
F(X) \otimes F(X) \otimes F(X) \oplus 3 F(X \vee X) \ar[d]^{\cong} \ar[r]^=
&
F(X) \otimes F(X) \otimes F(X) \oplus 3 F(X \vee X) \ar[d]^{\cong}\\
*\txt{$F(X \vee X) \otimes F(X) \oplus 4 F(X)$ \\
$\oplus$ 
$F(X \vee X)$} 
\ar[d]^{\cong}
&
*\txt{$F(X) \otimes F(X \vee X) \oplus 4 F(X)$\\ 
$\oplus$ 
$F(X \vee X)$}
\ar[d]^{\cong}\\
F(X \vee X \vee X) \oplus 3 F(X) \ar@{-->}[r]^p 
&
F(X \vee X \vee X) \oplus 3 F(X).
}
$$
The dotted arrow $p$ exists simply because all the other maps in the diagram 
are isomorphisms, and by description of the two vertical maps provided above, it is just  
the functor $F$ applied to a permutation of  
$X \vee X \vee X$. Of course, the maps $X \vee X \vee X \stackrel{+}{\rightarrow} X$
and $X \vee X \vee X \stackrel{perm}{\rightarrow} X \vee X \vee X 
\stackrel{+}{\rightarrow} X$  (where $perm$ is some permutation) are the same, and 
observing that the two maps are still equal after we apply $F$ to the maps 
(as $F$ is a functor), we have a commutative diagram:

$$
  \xymatrix{
F(X) \otimes F(X) \otimes F(X) \ar[d]
& \\
F(X) \otimes F(X) \otimes F(X) \oplus 3 F(X \vee X) \ar[d]^{\alpha \otimes 1} \ar[r]^=
&
F(X) \otimes F(X) \otimes F(X) \oplus 3 F(X \vee X) \ar[d]^{1 \otimes \alpha}\\
F(X \vee X \vee X) \oplus 3 F(X) \ar[r]^p \ar[dr]
&
F(X \vee X \vee X) \oplus 3 F(X) \ar[d]\\
&
F(X)
}
$$
which shows that the diagram \ref{expdiag} commutes and 
completes the proof of associativity of $\mu$. 
 
It remains to prove that the functor  $\tilde{F}$ from $\C$ to $\N$ that we just obtained
is coproduct preserving. 
Observe that it is immediate from the universal property of coproducts that 
$F(X) \otimes F(Y) \oplus F(X) \oplus F(Y)$ is the coproduct of $F(X)$ and $F(Y)$ in 
$\N$ and hence  $\tilde{F}$ is coproduct preserving.

\end{proof}

\section{Functors of Type $\frac{x}{1-x}$}

The next class of functors we consider mimics functions $f(x)=\frac{1}{1-ax}-1$. 
Here, as before, we subtract $1$ to make the functions reduced.
These functions are determined
by the functorial equation 
$$ f(x) +f(y) +2f(x)f(y) + f(x+y)f(x)f(y) =f(x+y) .$$
Moreover, every reduced analytic function that satisfies the above equation is of the 
form $f(x) = (\frac{1}{1-x}-1) \circ (ax)$.
This leads us to the following definition:

\begin{definition}
\label{def:ass}

A reduced functor $F$ from $\C$ to $Ch(\K)$ is of type $\frac{x}{1-x}$ if there 
is a natural isomorphism
\begin{eqnarray*}
\lefteqn{
\alpha_{X,Y}:F(X) \oplus F(Y) \oplus}\\ 
&&[F(X) \otimes F(Y)] \oplus [F(Y) \otimes F(X)]
\oplus [F(X \vee Y) \otimes F(X) \otimes F(Y)] \rightarrow   F(X \vee Y),
\end{eqnarray*}
such that the following natural diagram commutes:

$$
  \xymatrix{
*\txt{$F(X) \oplus F(Y) \oplus F(X) \otimes F(Y) \oplus F(Y) \otimes F(X)$ \\
 $\oplus$ \\ 
$F(X \vee Y) \otimes F(X) \otimes F(Y)$}
\ar[d]^{\cong} \ar[rr]^-{\alpha _{X,Y} }
&&
F(X \vee Y) \ar[d]^{\cong}\\
*\txt{$F(Y) \oplus F(X) \oplus F(Y) \otimes F(X) \oplus F(X) \otimes F(Y)$ \\
 $\oplus$ \\ 
$F(Y \vee X) \otimes F(Y) \otimes F(X)$} 
\ar[rr]^-{\alpha _{Y,X}} 
&&
F(Y \vee X)
}
$$
where the left vertical map takes the summand $F(X)$ to  $F(X)$, $F(Y)$ to  $F(Y)$,
$F(X) \otimes F(Y)$ and $F(Y) \otimes F(X)$ get mapped to $F(X) \otimes F(Y)$ and 
$F(Y) \otimes F(X)$ respectively,
and the map on the last component is the flip map $\tau$. 

\end{definition} 

Note that using the definition of cross effects (Section~\ref{sect:calc}) the functional equation
defining functors of type $ \frac{x}{1-x}$ can be rewritten as 
$$F(X) \otimes F(Y) \oplus F(Y) \otimes F(X)
  \oplus F(X \vee Y) \otimes F(X) \otimes F(Y) 
\stackrel{\cong}{\rightarrow}  cr_2 F(X,Y).$$
Similar restatements are true for all the other functors to be considered, and we will
refer to both equations as the defining equation of functors of the given type. 

\begin{example}

Let $\A$ be the category of non-unital associative differential graded $\K$-algebras.
(Clearly, $\A$ is naturally isomorphic to the category of $\T$ algebras in $Ch(\K)$, 
where $\T$ is the triple given by $\T = \bigoplus_{n \geq 1} X^{\otimes n}$ for all 
$X \in Ch(\K)$.) 
Here, and in future, the subscript $1$ refers to the fact that the summation starts 
at $n=1$ as opposed to $n=0$, and thus the functor is reduced.) 
Let $U_{\frac{x}{1-x}}$ be the forgetful functor 
$\A \rightarrow Ch(\K)$. We claim that $U_{\frac{x}{1-x}}$ is of type $\frac{x}{1-x}$. 
Observe, that if $X \vee Y$ is the coproduct of $X$ and $Y$ 
in $\A$, then so is $X \oplus Y \oplus X \otimes Y \oplus Y \otimes X \oplus (X \vee Y)
\otimes X \otimes Y$. Indeed, a pair of maps $X \rightarrow A$, $Y \rightarrow A$ in 
$\A$ induces a unique  map from each component of the sum above into $A$ in an obvious 
fashion. Now using the fact that two different coproducts of the same pair of objects 
must be isomorphic, we produce the necessary isomorphisms $\alpha_{X,Y}$. 

\end{example}

As in the case of exponential functors we have a lemma.

\begin{lemma}
\label{lemma:ass}

Every functor $F$ of type $\frac{x}{1-x}$ from $\C$ to $Ch(\K)$ can be factored 
as $F=U_{\frac{x}{1-x}} \circ A$, where $A$ is a reduced coproduct preserving functor 
from $\C$ to $\N$. 

\end{lemma}

\begin{proof}

Define $\mu : F(X) \otimes F(X) \rightarrow F(X)$ as the composite

\begin{eqnarray*}
\lefteqn{
F(X) \otimes F(X) \rightarrow
F(X) \oplus F(X) \oplus [F(X) \otimes F(X)] \oplus}\\ 
&&[F(X) \otimes F(X)]
  \oplus [F(X \vee X) \otimes F(X) \otimes F(X)] \rightarrow 
F(X \vee X) \rightarrow F(X)
\end{eqnarray*}

For the first map we have two choices and we simply fix one. The two choices of course
reflect the fact that any associative algebra $A$ admits two not necessarily isomorphic 
multiplications defining algebras $A$ and $A^{op}$. Naturally, the only thing we need to 
worry about is being consistent in our choices. 

As with exponential functors, we have a pair of commutative diagrams:

\noindent
$
  \xymatrix{
F(X \vee X) \otimes F(X) \ar[r]^-{F(+) \otimes 1 }\ar[d]^{\alpha_{X \vee X, X}}
&
F(X) \otimes F(X) \ar[d]^{\alpha_{X,X}} \\
F((X \vee X) \vee X) \ar[r]^-{F(+\vee 1)} &
F(X \vee X) 
}
$
\hspace{10pt}
$
  \xymatrix{
F(X) \otimes F(X \vee X) \ar[r]^-{1 \otimes F(+) }\ar[d]^-{\alpha_{X, X \vee X}}
&
F(X) \otimes F(X) \ar[d]^{\alpha_{X,X}} \\
F(X \vee (X \vee X)) \ar[r]^-{F(1\vee +)} &
F(X \vee X) 
}
$

We also have the following diagram:

$$
  \xymatrix{
F(X) \otimes F(X) \otimes F(X) \ar[d] \ar[r]^=
&
F(X) \otimes F(X) \otimes F(X) \ar[d]\\
F(X \vee X) \otimes F(X)  \ar[d]
&
F(X) \otimes F(X \vee X)  \ar[d]\\
F((X \vee X) \vee X)  \ar@{-->}[r]^p 
&
F(X \vee (X \vee X)) 
}.
$$

It is obtained just as the corresponding diagram in exponential case, except now
to get the two vertical maps we have to trace through the following (unfortunately 
very tedious) series of equivalences. (The point is, that, just like in exponential 
case, the left and right vertical composites differ essentially by a permutation 
of $X \vee X \vee X$.

$
 3 F( X \vee X) \oplus 4 F(X) \otimes F(X) \otimes F(X) \oplus
 2 F( X \vee X) \otimes F(X) \otimes F(X) \otimes F(X) \oplus
 F(X \vee X \vee X) \otimes   F( X \vee X) \otimes F(X) \cong 
$

$
 F(X \vee X) \oplus 4 F(X) \oplus 4 F(X) \otimes F(X) \oplus 
 2 F( X \vee X) \otimes F(X) \otimes F(X) 
 \oplus 4 F(X) \otimes F(X) \otimes F(X) \oplus
 2 F( X \vee X) \otimes F(X) \otimes F(x) \otimes F(X) \oplus
 F(X \vee X \vee X) \otimes   F( X \vee X) \otimes F(X) \cong
$

$
 F(X \vee X) \oplus 4 F(X) \oplus 2 F( X \vee X) \otimes F(X) \otimes F(X)
 \oplus 2 F(X) \otimes [2 F(X) \oplus 2 F(X) \otimes F(X) \oplus 
 F( X \vee X) \otimes F(X) \otimes F(X)] \oplus
 F(X \vee X \vee X) \otimes   F( X \vee X) \otimes F(X) \cong
$

$
 F(X \vee X) \oplus 4 F(X) \oplus 2 F( X \vee X) \otimes F(X) \otimes F(X)
 \oplus 2 F(X) \otimes F(X \vee X) \oplus 
 F(X \vee X \vee X) \otimes   F( X \vee X) \otimes F(X) \cong
$  
 
$
 3 F(X) \oplus 2 F( X \vee X) \otimes F(X) \otimes F(X) \oplus 
 F(X\vee X \vee X).
$

Again it is not hard to see that the map $p$ is the functor $F$ applied to a permutation 
of $X \vee X \vee X$, hence the two vertical arrows in the above diagram followed by 
$F(+)$ are equal, which proves that  $\mu$ is associative.

To see that the functor into the category of associative algebras that we just defined 
is coproduct preserving observe that $F(X \vee Y) \cong \colim C_i(X,Y)$, where the $C_i$'s 
are defined recursively:

$$C_0(X,Y) = F(X) \oplus F(Y) \oplus [F(X) \otimes F(Y)] \oplus 
                          [F(Y) \otimes F(X)]$$

$$C_{k+1}(X,Y) = C_k(X,Y) \oplus [C_k(X,Y) \otimes F(X) \otimes F(Y)].$$
Furthermore, clearly each pair of maps $f:F(X) \to A$ and  $g:F(Y) \to A$ in $\N$
induces a unique evident map $C_{k}(X,Y) \to A$ for all $k$. Thus, we obtain a (unique) 
map from  $F(X \vee Y) \cong \colim C_i(X,Y)$ to $A$, which proves by universal property 
of coproducts that  $F(X \vee Y)$ is the coproduct of $F(X)$ and $F(Y)$.

\end{proof}

\section{Taylor Towers of Functors of Type $\frac{x}{1-x}$}

In Calculus, Taylor series of functions can not always be computed. However, if the function 
is a solution to a sufficiently `nice' functional equation, its Taylor series can often be 
calculated. In this section, we take advantage of the fact that functors of type 
$\frac{x}{1-x}$ are defined via functional equations to compute their Taylor towers. 
We follow the approach taken by B. Johnson and R. McCarthy in ~\cite{Randy} in their treatment
of Exponential functors. First we will recall some of their results, needed for our 
computations. For completeness, their results on Taylor towers of exponential functors will 
also be presented here. 

We begin by observing that by Lemma 6.6 of ~\cite{Randy}, to understand the Taylor tower of
functors of type $\frac{x}{1-x}$ we only need to understand the Taylor tower of 
$U_ \frac{x}{1-x}$, since if $A$ is a reduced coproduct preserving functor, then 
$$P_n (F \circ A) \cong (P_n F) \circ A.$$
We list two results (Proposition 5.9 and Corollary 5.11 from ~\cite{Randy}) that are 
instrumental in our computations.

\begin{proposition}
\label{Randy1}
For a functor $F: \C \to Ch(\K)$, there is a natural isomorphism of $(n+1)$-multifunctors
$$\bigtriangledown ^n F(X_1, \cdots, X_n;Y) \cong
D_1 ^{(n)}[cr_{n+1}F(-, \cdots, -, Y) \oplus cr_n F(-, \cdots, -)](X_1, \cdots, X_n)$$
where the symbol  $D_1 ^{(n)}$ indicates that we are applying $D_1$ to each variable 
of the $n$-multifunctor separately, i.e. by holding all but one variable constant at 
all times.     
\end{proposition}

We do not present the definition of $\bigtriangledown ^n F$ in this document, 
as the description  of $\bigtriangledown ^n F$ provided by this Proposition  
is sufficient for our purposes. The functor $\bigtriangledown  F$ plays a role
analogous to that of the notion of a differential in classical Calculus. The following 
Corollary explains our interest in this functor.

\begin{corollary}
\label{Randy2}
For any functor $F: \C \to Ch(\K)$,
$$ D_n F \simeq \bigtriangledown ^n F(-, \cdots, -;0)_{h \Sigma_n}$$
\end{corollary}

Our objective is to describe the layers $D_n$ of the functor $U_\frac{x}{1-x}$
in terms of $D_1$. For the exponential functor $U_e$ this has been done in ~\cite{Randy},
where the authors computed the cross effects of $U_e$ and used the two results above to
conclude that 
\begin{equation}
\label{eq:ex}
D_n U_e \simeq \big( \bigotimes ^n D_1 U_e \big) _{h\Sigma_n}.
\end{equation}
Unfortunately, for the case of the functor $U_\frac{x}{1-x}$ computing the cross effects 
proved impractical. However, we are still able to get a formula analogous to the 
Equation~\ref{eq:ex}. 
First we demonstrate how to do it for $D_2$. 
By Proposition~\ref{Randy1}, 
\begin{eqnarray*}
\lefteqn{
  \bigtriangledown ^2 U_\frac{x}{1-x}(X_1, X_2;Y)}\\ 
   & & \cong 
   D_1 ^{(2)} [cr_3 U_\frac{x}{1-x}(-,-,Y) \oplus cr_2 U_\frac{x}{1-x}(-,-)](X_1, X_2) \\
   & & \cong 
   D_1 ^{(2)} [cr_3 U_\frac{x}{1-x}(-,-,Y)](X_1, X_2) 
\oplus D_1 ^{(2)} [cr_2 U_\frac{x}{1-x}(-,-)](X_1, X_2).
\end{eqnarray*} 
While computing $cr_3 U_\frac{x}{1-x}$ is not impossible, it is easily seen to produce 
long and tedious formulas which are hard to manage. Naturally, the higher cross effects 
which will be needed for the general case, are even worse. We make a few observations 
to simplify the computations. 

\begin{remark}
\label{rem:tayl1}
Recall that our intention is to apply Corollary~\ref{Randy2} to describe $D_2$. 
Consequently, we are interested in computing $D_1 ^{(2)} [cr_3 U_\frac{x}{1-x}(X_1,X_2,Y)]$
only for the special case $Y=0$. $cr_3 U_\frac{x}{1-x}(X_1,X_2,Y)$ is defined via the 
identity
$$cr_3 U_\frac{x}{1-x}(X_1,X_2,Y) \oplus cr_2 U_\frac{x}{1-x}(X_1,Y)
\oplus cr_2 U_\frac{x}{1-x}(X_2,Y) \cong cr_2 U_\frac{x}{1-x}(X_1 \vee X_2,Y).$$
It is immediate from the functional equation defining functors of type $\frac{x}{1-x}$
that the second and third terms on the left as well as the term on the right are 0 
when $Y=0$. Hence, for $Y=0$, $cr_3 U_\frac{x}{1-x}(X_1,X_2,Y) \cong 0$.
Thus, for $Y=0$, $D_1 ^{(2)} [cr_3 U_\frac{x}{1-x}(-,-,Y)](X_1, X_2) \cong 0$. 
 
\end{remark}

For our next observation, we need yet another result from ~\cite{Randy}. It is listed 
there as Lemma 3.3. 

\begin{lemma}
\label{Randy3}
Let $G: \C ^{\times n} \to Ch(\K)$ be an $n$-multireduced functor, i.e. $G$ is such that
$G(M_1, \cdots , M_n) \cong 0$ whenever any $M_i =0$. Then for 
$k<n$, $P_k G^{\bigtriangleup} \simeq 0$, where $G^{\bigtriangleup}$ is $G$ composed with 
the diagonal $\C \to \C ^{\times n}$. Consequently, for $k<n$, 
$D_k G^{\bigtriangleup} \simeq 0$.

\end{lemma}

\begin{remark}
\label{rem:tayl2}
Suppose $F: \C \to Ch(\K)$ is a reduced functor. Consider the functor
$$H = A \otimes F(X) \otimes B \otimes F(X) \otimes C,$$ 
where $A$, $B$, $C$ are objects 
in $Ch(\K)$. Then $H$ can be viewed as a functor $G:\C \times \C \to Ch(\K)$ composed 
with the diagonal. Since $F$ is reduced, $G$ is 2-multireduced. Hence by the above 
Lemma~\ref{Randy3}, $D_1 H = D_1 G^{\bigtriangleup} \simeq 0$.

\end{remark}

By Remark~\ref{rem:tayl1}, to understand $D_2 U_ \frac{x}{1-x}$, it is enough to 
consider 
$$D_1 ^{(2)}[cr_2 U_ \frac{x}{1-x} (-,-)](X_1, X_2).$$ 
Recall our formula for $ cr_2 U_ \frac{x}{1-x}$:
\begin{eqnarray}
cr_2 U_ \frac{x}{1-x} (X_1, X_2) 
& 
\cong 
&
U_ \frac{x}{1-x} (X_1) \otimes U_ \frac{x}{1-x} (X_2)
\oplus U_ \frac{x}{1-x} (X_2) \otimes U_ \frac{x}{1-x} (X_1) \nonumber\\
&& 
\oplus 
U_ \frac{x}{1-x} (X_1 \vee X_2) \otimes 
U_ \frac{x}{1-x} (X_1) \otimes U_ \frac{x}{1-x} (X_2) \nonumber \\
&
\cong
&
U_ \frac{x}{1-x} (X_1) \otimes U_ \frac{x}{1-x} (X_2)
\oplus U_ \frac{x}{1-x} (X_2) \otimes U_ \frac{x}{1-x} (X_1) \nonumber\\
&&
\oplus
[U_ \frac{x}{1-x} (X_1) \oplus U_ \frac{x}{1-x} (X_2) \nonumber \\
&&
\oplus
U_ \frac{x}{1-x} (X_1) \otimes U_ \frac{x}{1-x} (X_2) 
\oplus U_ \frac{x}{1-x} (X_2) \otimes U_ \frac{x}{1-x} (X_1) \nonumber\\
&&
\oplus 
U_ \frac{x}{1-x} (X_1 \vee X_2) \otimes 
U_ \frac{x}{1-x} (X_1) \otimes U_ \frac{x}{1-x} (X_2)] \nonumber \\
&&
\otimes U_ \frac{x}{1-x} (X_1) \otimes U_ \frac{x}{1-x} (X_2) \nonumber
\end{eqnarray}
Here to get the second equivalence we used the defining equation of functors of 
type $\frac{x}{1-x}$ again. 
Thus, all the terms except $U_ \frac{x}{1-x} (X_1) \otimes U_ \frac{x}{1-x} (X_2)$
and $U_ \frac{x}{1-x} (X_2) \otimes U_ \frac{x}{1-x} (X_1)$ have at least two 
factors of $U_ \frac{x}{1-x} (X_1)$ or two factors of $U_ \frac{x}{1-x} (X_2)$. 
Consequently by Remark~\ref{rem:tayl2}, these terms vanish after we apply either $D_1 ^1$ 
or $D_1 ^2$, where $D_1 ^i$ is $D_1$ applied to the $i'th$ variable with the other
variables held constant. Hence, 
\begin{eqnarray}
D_1 ^{(2)}[cr_2 U_ \frac{x}{1-x} (-,-)](X_1, X_2) 
& 
\cong 
&
D_1[U_ \frac{x}{1-x} (-_1) \otimes U_ \frac{x}{1-x} (-_2) \nonumber\\
&&
\oplus
U_ \frac{x}{1-x} (-_2) \otimes U_ \frac{x}{1-x} (-_1)](X_1, X_2) \nonumber \\
&
\cong 
&
D_1 U_ \frac{x}{1-x} (-)(X_1) \otimes D_1 U_ \frac{x}{1-x} (-)(X_2) \nonumber \\
&&
\oplus
D_1 U_ \frac{x}{1-x} (-)(X_2) \otimes D_1 U_ \frac{x}{1-x} (-)(X_1) \nonumber 
\end{eqnarray}
By Corollary~\ref{Randy2}, we conclude that
$$D_2 U_ \frac{x}{1-x} (-)(X) \cong 
D_1 U_ \frac{x}{1-x} (-)(X) \otimes D_1 U_ \frac{x}{1-x} (-)(X).$$
For the general case $D_n$, first we note that by induction
 $cr_{n+1} U_ \frac{x}{1-x}(X_1, \cdots, X_n, Y)$ can be written as a sum of terms
each of which contains at least one factor of $U_ \frac{x}{1-x}(Y)$. Hence, as in 
Remark~\ref{rem:tayl1}, for $Y=0$, 
$$D_1 ^{(n)}[cr_{n+1} U_ \frac{x}{1-x}(-, \cdots,-,Y)](X_1, \cdots, X_n) \simeq 0.$$ 
Thus, we have to worry only about  
$$D_1 ^{(n)}[cr_n U_ \frac{x}{1-x}(-, \cdots,-)](X_1, \cdots, X_n).$$ 
Observe that, analogous to the case $n=2$,  the only terms of 
$cr_n U_ \frac{x}{1-x}(X_1, \cdots,X_n)$ that do not 
contain at least two factors of $U_ \frac{x}{1-x}(X_i)$ for some $i$ are those of 
the form 
$U_ \frac{x}{1-x}(X_{\sigma (1)}) \otimes \cdots \otimes U_ \frac{x}{1-x}(X_{\sigma (n)})$, 
where $\sigma$ is any element of the permutation group $\Sigma _n$. 
Hence by Remark~\ref{rem:tayl2}, these are the only (possibly) non-vanishing terms after 
we apply $ D_1 ^{(n)}$. By Corollary~\ref{Randy2}, we are entitled to conclude that 
$$ D_n U_ \frac{x}{1-x} \cong 
 \bigotimes ^n D_1 U_ \frac{x}{1-x}.$$

\section{Logarithmic Functors and Their Taylor Towers}

From now on we will assume that $\K$ is a field of char 0.

In this section, we continue our investigation of functors patterned after elementary 
functions of one variable by considering logarithmic functions $-log(1-ax)$. While it is 
a good candidate for our purposes as it has a relatively simple Taylor series and is 
easily seen to be related to a common algebraic structure, we have yet another reason 
for singling out logarithmic functors.

To explain our interest 
we refer to the functions of real variables again. Observe that 
 $f(x)= \frac{1}{1-x}-1$ and $f(x)=e^x-1$  are respectively the first element and 
the limit of the sequence of functions 
$$f_n(x) = (1- \frac{x}{2^{n-1}})^{-2^{n-1}}-1.$$ 
Of course, here any number greater than 1 can be used instead of 2. 
This is reminiscent of the situation with the homology of little $n$-cubes 
operads $e_n$, where $e_n$ for $n>1$ interpolate between $e_1$ and $e_\infty$, which 
produce associative and commutative algebras respectively, much like the exponential 
functors and functors of type $\frac{x}{1-x}$. Hence the hope that the functors 
associated with the functions $f_n$ will provide an insight in the algebras over 
the homology of little $n$-cubes operads.     

However, before we attempt to pattern a class of functors after the functions $f_n$ above, 
we note that $f_n(x) = e^{-2^{n-1} \log(1- \frac{x}{2^{n-1}})}-1$. 
Thus, the function $-\log(1-x)$, in addition to being of interest on its own, is closely 
related to $f_n$. 

Before proceeding, we pause to observe 
that the suspension functor $\Sigma$ behaves like multiplication by $1/2$ (or any 
other number less than one) as it is linear and increases connectivity. In addition, 
multiplication by $2$ is represented by desuspension $\Sigma^{-1}$ since it is the 
inverse operation of $\Sigma$. 

As before, the following definition is inspired by the fact that $-\log(1-ax)$ are the 
only real variable analytic functions satisfying the functional equation

$$ f( \frac{x}{1-x}  \cdot \frac{y}{1-y}) + f(x) + f(y) = f(x + y).$$

\begin{definition}
\label{def:log}

A reduced functor $F:Ch(\K) \rightarrow Ch(\K)$ is logarithmic if there is a natural 
isomorphism

$$\alpha _{X,Y}: F(U X \otimes U Y) \oplus F(X) \oplus F(Y) \rightarrow F(X \vee Y)$$
where $U = \P_1 F$ and is a functor of type $\frac{x}{1-x}$,
such that for all $X, Y \in Ch(\K)$ the natural diagram commutes:

$$
  \xymatrix{
F(U X \otimes U Y) \oplus F(X) \oplus F(Y) \ar[r]^-{\alpha _{X,Y} }\ar[d]^\cong
&
F(X \vee Y) \ar[d]^{-1} \\
F(U Y \otimes U X) \oplus F(Y) \oplus F(X) \ar[r]^-{\alpha _{Y,X} } &
F(Y \vee X) 
},
$$
 
\end{definition} 

Before discussing the logarithmic functors, we digress a little to recall
the definition of a Lie algebra over chain complexes (i.e. in graded context).
The reason for introducing the Lie algebras at this point will be evident 
from the Example~\ref{ex:lie}.

\begin{definition}

$A \in Ch(\K)$ is a Lie algebra if it is equipped with a bracket operation
$[-,-]: A \otimes A \rightarrow A$, such that the composites 
$(A \otimes A)^{C_2} \rightarrow A \otimes A \rightarrow A$ and 
$(A \otimes A \otimes A)^{C_3} \rightarrow A \otimes A \otimes A \rightarrow A$  
are zero. Here $C_2$ and $C_3$ are the cyclic groups of size 2 and 3 respectively,
$A^G$ denotes the fixed points of the action of $G$ on $A$, and  the actions
of $C_2$ and $C_3$ are defined as follows:

The nontrivial element of $C_2$ acts on $A \otimes A$ by taking $a \otimes b$ to 
$(-1)^{|a||b|} b \otimes a$, and the two nontrivial elements of $C_3$ act by mapping 
$a \otimes b \otimes c$ to $(-1)^{|a||c|+|b||c|} c \otimes a \otimes b$ and 
$(-1)^{|a||c|+|a||b|} b \otimes c \otimes a$.

\end{definition}

\begin{example}
\label{ex:lie}

Let $\Lf$ be the free Lie algebra functor, i.e. the left adjoint to the forgetful 
functor $U_{\log}$ from the category $\L$ of Lie algebras to $Ch(\K)$. Alternatively,
$\Lf(M)$ can be defined as the Lie subalgebra of the Lie algebra $Lie({\mathbb T}(M))$
underlying the tensor algebra ${\mathbb T}(M)$, generated by $M$.  
We claim that $F=U_{\log} \circ \Lf$ is a logarithmic functor. 
First observe that, in characteristic $0$, by the Poincar\'e-Birkhoff-Witt Theorem, 
$\P_1 F \cong \T$, where $\P_1$ and $\T$ are the non-unital versions of 
the symmetric algebra triple $\P$ and the tensor algebra triple ${\mathbb T}$ respectively. 
Now it is not hard to see that the functor 
satisfies the equation
$$cr_2 U_{\log} \circ \Lf(X,Y) \cong U_{\log} \circ \Lf (\T X \otimes \T Y).$$

\end{example}

As a curious point, it is worth noting that the Poincar\'e-Birkhoff-Witt Theorem
is a restatement of the simple functional identity
$$e ^{-\log(1-x)} -1 = \frac{x}{1-x}.$$
As with functors of type $\frac{x}{1-x}$, our next goal is to express the higher
derivatives of logarithmic functors in terms of $D_1$. Again, this is an application
of the theory developed by B.Johnson and R.McCarthy in ~\cite{Randy}. 

Let $F$ be a logarithmic functor. We begin with the computation of $D_2 F$. 
By Proposition~\ref{Randy1}, 
\begin{eqnarray}
  \bigtriangledown ^2 F(X_1, X_2;Y) & \cong &
   D_1 ^{(2)} [cr_3 F(-,-,Y) \oplus cr_2 F(-,-)](X_1, X_2)
  \nonumber \\
   &
   \cong &
   D_1 ^{(2)} [cr_3 F(-,-,Y)](X_1, X_2) 
\oplus D_1 ^{(2)} [cr_2 F(-,-)](X_1, X_2) \nonumber .
\end{eqnarray}
Recall that 
$$cr_3 F(X_1, X_2, Y) \oplus cr_2 F(X_1, Y) \oplus cr_2 F(X_2, Y) 
\cong cr_2 F(X_1 \vee X_2, Y).$$
However, $ D_1 ^{(2)} cr_2 F(X_i, Y) =0$ for $i=1,2$, since $cr_2 F(X_i, Y) =0$ is 
constant with respect to $X_1$ or $X_2$. Hence,
\begin{eqnarray}
D_1 ^{(2)} [cr_3 F(-,-,Y)](X_1, X_2)  
& \cong & 
D_1 ^{(2)} [cr_2 F(- \vee -, Y)](X_1, X_2)
\nonumber\\
& \cong & 
D_1 ^{(2)} [F(U(- \vee -) \otimes UY)](X_1, X_2).
\nonumber
\end{eqnarray}
Of course, our intention is to apply the Corollary~\ref{Randy2}, and hence we are interested 
only in the case $Y=0$. Observe that since the functors $A \otimes -$, $U$ and $F$ are reduced,
$D_1 ^{(2)} [F(U(- \vee -) \otimes UY)](X_1, X_2) \cong 0$ for $Y=0$. Thus, we only need to 
compute $D_1 ^{(2)} [cr_2 F(-,-)](X_1, X_2)$.

To do this we use the result listed as Lemma 5.7 in ~\cite{Randy} which establishes a 
chain rule for functors analogous to the chain rule for functions of one variable. 

\begin{lemma}
\label{chain}
Let ${\mathcal A}$ and ${\mathcal A}^{'}$ be Abelian categories, 
and let functors $F$ and $G$ be such that $G: \C \to  {\mathcal A}^{'}$ and 
$F:{\mathcal A}^{'} \to Ch{\mathcal A}$. Then
$$D_1(F \circ G) \cong D_1 F \circ D_1 G.$$
\end{lemma} 

We use this lemma to compute the differential $D_1 ^{(2)}$ of $cr_2 F$:
\begin{eqnarray}
D_1 cr_2F(-,X_2) (X_1)
& \cong & 
D_1 F(U(-) \otimes UX_2)(X_1) \cong D_1 [F \circ (U(-) \otimes UX_2)](X_1)
\nonumber\\
&\cong &
D_1 F \circ D_1[U(-) \otimes UX_2](X_1) 
\nonumber\\
& \cong & 
D_1 F \circ [D_1 U(-)(X_1) \otimes UX_2]
\nonumber\\
& \cong &
[D_1 U(-)(X_1) \otimes UX_2] \otimes D_1 F(\K),
\nonumber
\end{eqnarray}
where $\K$ is the ground field. The last equivalence holds because $D_1 F$ is linear. 
Furthermore, note that
\begin{equation}
\label{eq:du}
D_1 U(-)(X_1) \cong D_1(\P_1 F) U(-)(X_1) \cong D_1 \P_1 \circ D_1 F(-)(X_1)
\cong D_1 F(-)(X_1).
\end{equation}
Combining these calculations we have that
$$D_1 cr_2F(-,X_2) (X_1) \cong D_1 F(X_1) \otimes U(X_2) \otimes D_1 F(\K).$$
Consequently, using the linearity of $D_1$ again, we get
\begin{eqnarray}
D_1 ^{(2)} cr_2F(-,-) (X_1, X_2) 
& \cong &
D_1 ^2 [D_1 F(X_1) \otimes U(-) \otimes D_1 F(\K)](X_2)
\nonumber \\
& \cong &
D_1 F(X_1) \otimes D_1 U(-)(X_2) \otimes D_1 F(\K)
\nonumber \\
& \cong & 
D_1 F(X_1) \otimes D_1 F(X_2) \otimes D_1 F(\K).
\nonumber
\end{eqnarray}
Thus, we have completed the computation of $\bigtriangledown ^2 F(X,X;0)$:

$$\bigtriangledown ^2 F(X,X;0) \cong D_1 F(X_1) \otimes D_1 F(X_2) \otimes D_1 F(\K).$$

Now we consider the general case. 

\begin{remark}
\label{rem:cross}
Recall that 
\begin{eqnarray}
cr_{n+1} F(X_1, X_2, \cdots , X_n, Y) 
& \oplus & 
cr_n F(X_2, \cdots , X_n, Y) 
\oplus 
cr_n F(X_1, X_3, \cdots , X_n, Y) 
\nonumber \\
& \cong & 
cr_n F(X_1 \vee X_2, X_3, \cdots , X_n, Y)
\nonumber
\end{eqnarray}
Hence when we apply $D_1 ^{(n)}$ to the two sides of the above equation, the second 
and the third terms on the left vanish since they are constant with respect to 
$X_1$ and $X_2$ respectively. Thus,
$$ D_1 ^{(n)} cr_{n+1} F(-, \cdots, -, Y) (X_1, \cdots , X_n) \cong
D_1 ^{(n)} cr_{n} F(-\vee -,-, \cdots, -, Y) (X_1, \cdots , X_n).$$
\end{remark}
As with functors of type $\frac{x}{1-x}$, we are interested in the case $Y=0$ (see 
Corollary~\ref{Randy2}). By the above arguments,
$$D_1 ^{(n)} cr_{n} F(-\vee -,-, \cdots, -, Y) \cong
D_1 ^{(n)} cr_{n-1} F(-\vee -\vee -, \cdots, -, Y).$$
Hence, for $Y=0$, by induction on $n$, we have that 
$$D_1 ^{(n)} cr_{n} F(-\vee -,-, \cdots, -, Y) (X_1, \cdots , X_n) \cong 0.$$
Consequently, for $Y=0$
$$D_1 ^{(n)} cr_{n+1} F(-, \cdots, -, Y) (X_1, \cdots , X_n) \cong 0.$$
Thus, we only need to understand $D_1 ^{(n)} cr_{n} F(-, \cdots, -)(X_1, \cdots , X_n)$. 
By iterating the argument in Remark~\ref{rem:cross}, we get that 
\begin{eqnarray}
D_1 ^{(n)} cr_{n} F(-, \cdots, -)(X_1, \cdots , X_n) 
&\cong&
D_1 ^{(n)} cr_2 F(- \vee \cdots \vee -,-)(X_1, \cdots , X_n)
\nonumber\\
&\cong&
D_1 ^{(n)} F(U(- \vee \cdots \vee -) \otimes U(-))(X_1, \cdots , X_n)
\nonumber
\end{eqnarray}
Recall that $U$ is a functor of type $\frac{x}{1-x}$, and thus using the computations 
performed in the previous section, by Equation~\ref{eq:du}, we have that
\begin{eqnarray*}
\lefteqn{
D_1 ^{(n-1)} U(- \vee \cdots \vee -)(X_1, \cdots , X_{n-1}) } \\
& & \simeq  
 \bigoplus   
D_1 U(-)(X_{\sigma(1)}) \otimes \cdots \otimes D_1 U(-)(X_{\sigma(n-1)}) \\
& & \simeq 
 \bigoplus  
D_1 F(-)(X_{\sigma(1)}) \otimes \cdots \otimes D_1 F(-)(X_{\sigma(n-1)}) 
\end{eqnarray*} 
where the sums on the right are taken over all permutations $\sigma$ 
of the group $\Sigma_{n-1}$.
Combining these calculations with the Lemma~\ref{chain} on chain rule, as well as 
our computation of the second layer $D_2$, we get that
\begin{eqnarray*}
\lefteqn{
D_1 ^{(n)} cr_n F(X_1, \cdots , X_n)} \\
&& \simeq 
D_1 ^n[D_1 ^{(n-1)} cr_n F(-, \cdots ,-,-)(X_1, \cdots , X_{n-1})](X_n) \\
&& \simeq 
D_1 ^n[D_1 ^{(n-1)} F(U(- \vee \cdots \vee -) \otimes U(-))(X_1, \cdots , X_{n-1})](X_n) \\
&& \simeq 
[\bigoplus
D_1 F(-)(X_{\sigma(1)}) \otimes \cdots \otimes D_1 F(-)(X_{\sigma(n-1)})] 
\otimes
D_1 ^n [D_1 F (\K \otimes  U(-))](X_n) \\
&& \simeq 
[\bigoplus
D_1 F(-)(X_{\sigma(1)}) \otimes \cdots \otimes D_1 F(-)(X_{\sigma(n-1)})]  
\otimes
D_1 F(X_n) \otimes D_1 F(\K) 
\end{eqnarray*} 
Hence, we are allowed to conclude that
$$\bigtriangledown ^n F(X, \cdots ,X; 0) \cong 
\bigoplus [D_1 F(X) \otimes \cdots \otimes D_1 F(X)] \otimes D_1 F(X) \otimes D_1 F(\K) $$
where there are $n-1$ factors in the brackets of the sum on the right, and the sum
is again taken over all permutations of the group $\Sigma_{n-1}$. 

Hence for $n>1$, by Corollary~\ref{Randy2}, $D_n F$ is the tensor product over 
$\Sigma_n$ of the $n'th$ tensor power of $D_1 F$ with an $n-1$ dimensional 
representation of $\Sigma_n$. Combining this with the nature of the action of 
$\Sigma_n$ on powers of $D_1 F$ described above, we conclude that

$$D_n F \cong Lie(n) \otimes_{h \Sigma_n} \bigotimes^n D_1 FU \otimes D_1 F(\K),$$ 
where $Lie(n)$ is the Lie representation of $\Sigma_n$.

\section{Functors of Type $f_n$}

Now we get back to  functors of type $f_n$ for $n>1$. The functional equation (for
reduced analytic functions) defining $f_n(ax)$ is as follows:

$$ f(x+y) = f(x) + f(y) + f(x) f(y) + 
f \big ( 2^{n-1} \cdot \frac{x/2^{n-1}}{1-x/2^{n-1}} \cdot \frac{y/2^{n-1}}{1-y/2^{n-1}} \big)$$ 

$$+ f \big( 2^{n-1} \cdot \frac{x/2^{n-1}}{1-x/2^{n-1}} \cdot \frac{y/2^{n-1}}{1-y/2^{n-1}} \big) 
\cdot [ f(x) + f(y) + f(x) f(y)].$$

Recalling that multiplication by $2$ corresponds to desuspension, we give the following
definition:  

\begin{definition}
\label{def:f_n}

A reduced functor $F : Ch(\K) \rightarrow Ch(\K)$ is of type $f_n$,  if 
there is a natural isomorphism 
\begin{eqnarray*}
\lefteqn{
\alpha_{X,Y}: F(X) \oplus F(Y)  \oplus F(X) \otimes F(Y) \oplus} \\
&& F(\Sigma^{-n+1} U X \otimes U Y) \oplus \\
&& F(\Sigma^{-n+1} U X \otimes U Y) \otimes
[F(X) \oplus F(Y)  \oplus F(X) \otimes F(Y)] \rightarrow
F(X \vee Y)
\end{eqnarray*}
where $U=\P_1 \Sigma^{n-1}F$ is a functor of type $\frac{x}{1-x}$, 
such that the following two diagrams commute:
 
$$ \xymatrix{
F(X) \otimes F(Y)  \ar[r]^-{\alpha _{X,Y} }\ar[d]^\cong
&
F(X \vee Y) \ar[d]^\cong \\
F(Y) \otimes F(X)  \ar[r]^-{\alpha _{Y,X} } &
F(Y \vee X) 
}$$ 
$$ \xymatrix{
F(\Sigma^{-n+1} U X \otimes U Y) \ar[r]^-{\alpha _{X,Y} }\ar[d]^\cong
&
F(X \vee Y) \ar[d]^{-1} \\
F(\Sigma^{-n+1} U Y \otimes U X) \ar[r]^-{\alpha _{Y,X} } &
F(Y \vee X). 
}$$
\end{definition}
Note that the ingredients of the functional equation defining the functors of type 
$f_n$ are familiar. In particular, the term $F(\Sigma^{-n+1} U X \otimes U Y)$ 
resembles the defining equation of logarithmic functors with a shift, while 
the term $F(X) \otimes F(Y)$ is the key component of exponential functors.
Hence it is reasonable to expect that examples of functors of type $f_n$ evaluated 
at an object $X$, are equipped with both a Lie bracket and a commutative multiplication. 

We make a few observations before we can produce examples of functors of type $f_n$. 

\begin{remark}
\label{rem:f_n}
The purpose of this remark is to observe that if $L$ is a logarithmic functor, then 
$\P \Sigma^{-n+1} L \Sigma^{n-1} X$ is equipped with a Lie bracket of degree $n-1$,
in particular, we produce a map:
\begin{equation}
\label{eq:bracket}
[-,-]_{n-1}: \P (\Sigma^{1-n} L (\Sigma^{n-1} X)) \otimes \P (\Sigma^{1-n} L (\Sigma^{n-1} X))
 \rightarrow \Sigma^{1-n} \P (\Sigma^{1-n} L (\Sigma^{n-1} X)).
\end{equation}
A typical summand on the left is of the form 
\begin{equation}
\label{eq:summand} 
[\Sigma^{1-n} L (\Sigma^{n-1} X)]^{\otimes l}/ \Sigma_l \otimes 
 [\Sigma^{1-n} L (\Sigma^{n-1} X)]^{\otimes k}/ \Sigma_k .
\end{equation} 
So to describe (~\ref{eq:bracket}),
it is enough to produce the maps
\begin{eqnarray}
\label{eq:brac}
\lefteqn{
 [\Sigma^{1-n} L (\Sigma^{n-1} X)]^{\otimes l}/ \Sigma_l \otimes 
  [\Sigma^{1-n} L (\Sigma^{n-1} X)]^{\otimes k}/ \Sigma_k} \\
&& \to 
\Sigma^{1-n} [\Sigma^{1-n} L (\Sigma^{n-1} X)]^{\otimes k+l-1}/ \Sigma_{k+l-1}.
\nonumber
\end{eqnarray}
First, since $L$ is logarithmic,  we have the bracket 
$[-,-]_L: L (\Sigma^{n-1} X) \otimes L (\Sigma^{n-1} X) \rightarrow L (\Sigma^{n-1} X)$, 
which (by desuspending $2(n-1)$ times) produces a map 
$ \Sigma^{1-n} L (\Sigma^{n-1} X) \otimes \Sigma^{1-n} L (\Sigma^{n-1} X) 
\rightarrow \Sigma^{1-n} (\Sigma^{1-n} L (\Sigma^{n-1} X))$. As stated, our objective is to 
extend this map to terms described by Equation~\ref{eq:summand}.
This will be done in two steps. First we consider the case $l=1$. To begin, we make a 
general observation. 

Fix a chain complex $C$. Consider the map
$$C^{\otimes k} \to \bigoplus _1 ^k [C \otimes C^{\otimes (k-1)}/\Sigma_{k-1}] 
\stackrel{+}{\to} C \otimes C^{\otimes (k-1)}/\Sigma_{k-1},$$
where the second map is the obvious addition map, and the first map is defined on 
an element $a_1 \otimes \cdots \otimes a_k \in C^{\otimes k}$ by 
$$a_1 \otimes \cdots \otimes a_k  \mapsto 
\Sigma_\sigma a_{\sigma(1)} \otimes \cdots \otimes a_{\sigma (k)},$$
with the sum on the right ranging over $k$ cyclic permutations of $k$ letters.  
Note that this map factors through $C^{\otimes k}/\Sigma_{k}$, thus producing a 
map $\gamma _C :C^{\otimes k}/\Sigma_{k} \to C \otimes C^{\otimes (k-1)}/\Sigma_{k-1}$.

Returning to our case, observe that we have a sequence of maps:

\begin{equation}
\label{eq:brac2}
\begin{array}{llll}
\Sigma^{1-n} L (\Sigma^{n-1} X) \otimes 
[\Sigma^{1-n} L (\Sigma^{n-1} X)]^{\otimes k} /\Sigma_{k} \\
\rightarrow
\Sigma^{1-n} L (\Sigma^{n-1} X) \otimes \Sigma^{1-n} L (\Sigma^{n-1} X) \otimes
[\Sigma^{1-n} L (\Sigma^{n-1} X)]^{\otimes k-1}/\Sigma_{k-1} \\
\stackrel{[-,-] \otimes 1}{\longrightarrow}
\Sigma^{1-n} (\Sigma^{1-n} L (\Sigma^{n-1} X)) \otimes
[\Sigma^{1-n} L (\Sigma^{n-1} X)]^{\otimes k-1}/\Sigma_{k-1}\\ 
\rightarrow
\Sigma^{1-n} [\Sigma^{1-n} L (\Sigma^{n-1} X)]^{\otimes k}/ \Sigma_k,
\end{array}
\end{equation}
where the first map is the map $\gamma_{\Sigma^{1-n} L (\Sigma^{n-1} X)}$ defined above.
The general case follows by induction:
$$
\begin{array}{llll}
[\Sigma^{1-n} L (\Sigma^{n-1} X)]^{\otimes l}/ \Sigma_l \otimes 
[\Sigma^{1-n} L (\Sigma^{n-1} X)]^{\otimes k}/ \Sigma_k \\
\stackrel{\gamma \otimes 1} {\rightarrow}
[\Sigma^{1-n} L (\Sigma^{n-1} X)]^{\otimes l-1}/ \Sigma_{l-1} \otimes 
\Sigma^{1-n} L (\Sigma^{n-1} X) \otimes
[\Sigma^{1-n} L (\Sigma^{n-1} X)]^{\otimes k}/ \Sigma_k \\
\stackrel{1 \otimes [-,-]}{\longrightarrow}
[(\Sigma^{1-n} L (\Sigma^{n-1} X))]^{\otimes l-1}/ \Sigma_{l-1} \otimes
\Sigma^{1-n}[\Sigma^{1-n} L (\Sigma^{n-1} X)]^{\otimes k}/ \Sigma_k \\
\rightarrow
\Sigma^{1-n} [\Sigma^{1-n} L (\Sigma^{n-1} X)]^{\otimes k+l-1}/ \Sigma_{k+l-1} 
\end{array}
$$
The map $[-,-]_{n-1}$ we just defined is a Lie bracket as it is essentially
the $n-1$ suspension of the Lie bracket  $[-,-]_L$.

\end{remark}

This prompts us to recall the definition of $n$-Poisson algebras
(see ~\cite{Getz}).

\begin{definition}

$P \in Ch(\K)$ is an $n- Poisson$ algebra if it is equipped with a commutative 
multiplication $\mu: P \otimes P \rightarrow P$ and a bracket 
$[-,-]: P \otimes P \rightarrow  \Sigma^{-n} P$, such that, as with Lie algebras,
the composites $(P \otimes P)^{C_2} \rightarrow  P \otimes P \rightarrow  \Sigma^{-n}P$
and $(P \otimes P \otimes P)^{C_3} \rightarrow  P \otimes P \otimes P
\rightarrow  P \otimes \Sigma^{-n}P \rightarrow  \Sigma^{-2n} P$ are $0$ and the 
following diagram commutes:

$$
  \xymatrix{
P \otimes (P \otimes P) \ar[r]^{1 \otimes N} \ar[dd]^{1 \otimes \mu}
&
P \otimes (P \otimes P)^{C_2} \ar[r]
&
P \otimes (P \otimes P) \ar[r]
&
(P \otimes P) \otimes P \ar[d]^{[-,-] \otimes 1}\\
&&& 
\Sigma^{-n} P \otimes P \ar[d]^{\Sigma^{-n} \mu}\\
P \otimes P \ar[rrr]^-{[-,-]}
&
&
&
\Sigma^{-n} P 
}
$$
where $N$ is the Norm map of the action by $C_2$. The group actions are as before.  
Clearly, this diagram produces the usual Poisson relation. 

\end{definition}

\begin{example}
\label{example:f_n}
Here we produce examples of functors of type $f_n$. 
As it was mentioned earlier, the real valued exponential and logarithmic functions
are related to the functions $f_n$ via the formula 
$f_n(x) = e^{-2^{n-1} \log(1- \frac{x}{2^{n-1}})}-1$. 
The functorial analogue of the right hand side is, of course, 
$F_n(X):= E (\Sigma^{1-n} L (\Sigma^{n-1} X))$, where $E$ is an exponential functor
and $L$ logarithmic. Applying Definitions~\ref{def:log} and ~\ref{def:exp}
successively we get a series of isomorphisms:
\begin{eqnarray*}
\lefteqn{
F_n(X \vee Y)= E [\Sigma^{1-n} L (\Sigma^{n-1} (X \vee Y))]} \\
&& \cong 
E [\Sigma^{1-n} (L (\Sigma^{n-1} X) \oplus L(\Sigma^{n-1} Y) \oplus 
L(U_L \Sigma^{n-1} X \otimes U_L \Sigma^{n-1} Y))] \\
&&
\cong
E [\Sigma^{1-n} L (\Sigma^{n-1} X)] \oplus E [\Sigma^{1-n} L (\Sigma^{n-1} Y)]\\ 
&&
\oplus
E [\Sigma^{1-n} L(U_L \Sigma^{n-1} X \otimes U_L \Sigma^{n-1} Y)] \oplus
E [\Sigma^{1-n} L (\Sigma^{n-1} X)] \otimes E [\Sigma^{1-n} L (\Sigma^{n-1} Y)] \\
&&
\oplus
E [\Sigma^{1-n} L (\Sigma^{n-1} X)] \otimes
E [\Sigma^{1-n} L(U_L \Sigma^{n-1} X \otimes U_L\Sigma^{n-1}  Y)] \\
&&
\oplus
E [\Sigma^{1-n} L (\Sigma^{n-1} Y)] \otimes
E [\Sigma^{1-n} L(U_L\Sigma^{n-1}  X \otimes U_L\Sigma^{n-1}  Y)] \\
&&
\oplus
E [\Sigma^{1-n} L (\Sigma^{n-1} X)] \otimes    
E [\Sigma^{1-n} L (\Sigma^{n-1} Y)] \otimes
E [\Sigma^{1-n} L(U_L \Sigma^{n-1} X \otimes U_L \Sigma^{n-1} Y)] \\
&&
\cong
F_n(X) \oplus F_n(Y) \oplus F_n( \Sigma^{1-n} U_L \Sigma^{n-1} X \otimes U_L \Sigma^{n-1} Y)
\oplus
F_n(X) \otimes F_n(Y) \\
&&
\oplus
F_n(X) \otimes 
F_n( \Sigma^{1-n} U_L \Sigma^{n-1} X \otimes U_L \Sigma^{n-1} Y) \\
&&\oplus
F_n(Y) \otimes F_n(\Sigma^{1-n} U_L \Sigma^{n-1} X \otimes U_L \Sigma^{n-1} Y) \\
&&
\oplus
F_n(X) \otimes F_n(Y) \otimes F_n(\Sigma^{1-n} U_L \Sigma^{n-1} X \otimes U_L \Sigma^{n-1} Y).
\end{eqnarray*}

\noindent
Here by $U_L$ we denote the functor of type $\frac{x}{1-x}$ associated with the logarithmic functor 
$L$. In section on logarithmic functors, this functor is denoted simply by $U$. Here and in 
what follows we use the subscripts to avoid confusion.

Now we specialize to the case $E=\P$ and $L=\Lf$. To see that 
$F_n$ satisfies the defining equation of functors of type $f_n$ for this case, we need to 
show that $U_{F_n}X \cong U_L \Sigma^{n-1}X$. Observe that 
$U_{F_n}X \cong \P_1 \Sigma^{n-1} \P \Sigma^{1-n} \Lf \Sigma^{n-1}X$ and 
$U_L \Sigma^{n-1}X \cong \P_1 \Lf (\Sigma^{n-1}X)$, hence the desired identity follows
from the Poincar\'e-Birkhoff-Witt theorem and the
fact that the Lie algebras $\Sigma^{n-1} \P \Sigma^{1-n} \Lf \Sigma^{n-1}X$ and 
$\Lf (\Sigma^{n-1}X)$ have isomorphic universal enveloping algebras (see also 
Propositions 3.8 and 3.9 of~\cite{Milnor}). 

We note that these functors $F_n$ take values in the category of $n$-Poisson algebras.
An explicit description of the bracket multiplication $[-,-]_{n-1}$ of $F_n$ was given 
in Remark~\ref{rem:f_n}, as for the Poisson condition,
it is evident from the process which we used to extend the map 
$ \Sigma^{1-n} L (\Sigma^{n-1} X) \otimes \Sigma^{1-n} L (\Sigma^{n-1} X) 
\rightarrow \Sigma^{1-n} (\Sigma^{1-n} L (\Sigma^{n-1} X))$ to  
$\P (\Sigma^{1-n} L (\Sigma^{n-1} X)) \otimes \P (\Sigma^{1-n} L (\Sigma^{n-1} X))$.

\end{example}

\section{Taylor Towers of Functors of Type $f_n$}

Now we turn to computing the Taylor tower of the functor $F$ of type $f_n$. As before, we are 
going to use induction to compute the layers $D_k$. Hence we begin with $D_2$. 
The approach is the same: we start by analyzing the differential 
$\bigtriangledown ^2 F(X_1, X_2;Y)$. Recall that by Proposition~\ref{Randy1},
$$\bigtriangledown ^2 F(X_1, X_2;Y) \cong  
D_1 ^{(2)} [cr_3 F (-,-,Y)](X_1, X_2) 
\oplus D_1 ^{(2)} [cr_2 F(-,-)](X_1, X_2).$$
By the argument employed when discussing the logarithmic functors, we have that
\begin{eqnarray*}
\lefteqn{
D_1 ^{(2)} [cr_3 F (-,-,Y)](X_1, X_2) \cong   
D_1 ^{(2)} [cr_2 F (-\vee-,Y)](X_1, X_2)} \\
&& \cong 
D_1 ^{(2)}[F (-\vee-)\otimes F(Y) \oplus F(\Sigma^{-n+1} U(-\vee-)\otimes U(Y))  \\
&& \oplus
F(\Sigma^{-n+1} U(-\vee-)\otimes U(Y))  \\
&& \otimes 
(F (-\vee-) \oplus F(Y) \oplus F(-\vee-)\otimes F(Y))](X_1, X_2).
\end{eqnarray*}
Recall that we are interested in this object only for the case $Y \cong 0$ 
(see Corollary~\ref{Randy2}). Observe that 
$$ D_1 ^{(2)}[F (-\vee-)\otimes F(Y)] (X_1, X_2)
\cong  D_1 ^{(2)}[F (-\vee-)](X_1, X_2) \otimes F(Y) \cong 0$$
for $Y \cong 0$. The term  
$D_1 ^{(2)} [F(\Sigma^{-n+1} U(-\vee-)\otimes U(Y))](X_1, X_2)$ also vanishes 
when evaluated at $Y \cong 0$ by chain rule (Lemma~\ref{chain}), and since 
the functors $F$, $U$, $\Sigma$ and $A \otimes -$ (for all $A$) are reduced. 
Finally, the term 
$$D_1 ^{(2)}F(\Sigma^{-n+1} U(-\vee-)\otimes U(Y)) \otimes
(F (-\vee-) \oplus F(Y) \oplus F(-\vee-)\otimes F(Y))](X_1, X_2)$$
is also equivalent to $0$, by chain rule (Lemma~\ref{chain}). 
Thus, for $Y \cong 0$,
$$D_1 ^{(2)} [cr_3 F (-,-,Y)](X_1, X_2) \cong 0.$$
Hence it remains to compute $D_1 ^{(2)} [cr_2 F(-,-)](X_1, X_2)$. 

Observe that each summand of the term 
$$F(\Sigma^{-n+1} U X_1 \otimes U X_2) \otimes
[F(X_1) \oplus F(X_2)  \oplus F(X_1) \otimes F(X_2)] \rightarrow
F(X \vee Y)$$ 
is 2-reduced in either $X_1$, $X_2$ or both. Consequently, $D_1 ^{(2)}$ applied to this 
term is equivalent to $0$. Hence,
$$D_1 ^{(2)} [cr_2 F(-,-)](X_1, X_2) \cong
D_1 ^{(2)}[F(-) \otimes F(-) \oplus F(\Sigma^{-n+1} U(-)\otimes U(-))](X_1, X_2).$$
First, observe that 
$$D_1 U(-)(X) \cong D_1 (\P_1 \Sigma ^{n-1} F)(-)(X) \cong \Sigma ^{n-1} D_1 F(X).$$
The following computations are analogous to those performed for logarithmic functors,
so we omit explanations and refer to the logarithmic case for details:
\begin{eqnarray}
D_1 ^{(2)} F(\Sigma^{-n+1} U(X_1)\otimes U(X_2)) 
& \cong &
D_1 ^2 [D_1 F(\Sigma^{-n+1} D_1 U(-)(X_1) \otimes U(-))](X_2) \nonumber \\
& \cong &
\Sigma^{-n+1} D_1 U(X_1) \otimes D_1 ^2[ D_1 F(\K \otimes U(-))](X_2) \nonumber \\
& \cong &
\Sigma^{-n+1} \Sigma ^{n-1} D_1 F(X_1) \otimes D_1 ^2[ D_1 F(\K \otimes U(-))](X_2)
\nonumber \\
& \cong &
D_1 F(X_1) \otimes D_1 U(-)(X_2) \otimes D_1 F(\K) \nonumber \\
& \cong &
\Sigma ^{n-1} D_1 F(X_1) \otimes D_1 F(X_2) \otimes D_1 F(\K) \nonumber 
\end{eqnarray}  
Using these computations along with the Equation~\ref{eq:ex} of Section 6 which  
describes the layers of exponential functors, 
we are allowed to conclude that 
$$\bigtriangledown ^2 F(X, X;0) \cong 
D_1 F(X) \otimes D_1 F(X) \oplus 
\Sigma ^{n-1} D_1 F(X) \otimes D_1 F(X) \otimes D_1 F(\K)$$

For the general case of the $k'th$ layer, observe that by Remark~\ref{rem:cross}
\begin{eqnarray*}
\lefteqn{
D_1 ^{(k)} cr_{k+1} F(-, \cdots, -, Y) (X_1, \cdots , X_k)} \\ 
&& \cong 
D_1 ^{(k)} cr_{k} F(-\vee -,-, \cdots, -, Y) (X_1, \cdots , X_k)  \\
&& \cong 
\cdots \cong D_1 ^{(k)} cr_{2} F(-\vee \cdots \vee -, Y) (X_1, \cdots , X_k).
\end{eqnarray*}

Hence for $Y \cong 0$, by the above argument for the second differential 
$\bigtriangledown ^2 F$, we have that 
$$D_1 ^{(k)} cr_{2} F(-\vee \cdots \vee -, Y) (X_1, \cdots , X_k) \cong 0.$$
Consequently, for $Y \cong 0$, 
$$D_1 ^{(k)} cr_{k+1} F(-, \cdots, -, Y) (X_1, \cdots , X_k) \cong 0.$$
Thus, by Proposition~\ref{Randy1},
\begin{eqnarray}
\bigtriangledown ^k F(X_1, \cdots , X_k; 0) 
& \cong & 
D_1 ^{(k)} [cr_{k} F(-, \cdots, -)](X_1, \cdots , X_k) \nonumber \\
& \cong &
D_1 ^{(k)}[cr_2 F(- \vee \cdots \vee -, -)](X_1, \cdots ,X_{k-1}, X_k),\nonumber
\end{eqnarray}
where the second equivalence is by Remark~\ref{rem:cross}. To ease the notation, we define 
$A_k (X_1, \cdots , X_k)$ for $k>1$ to be 
\begin{equation}
\label{eq:A_k}
A_k (X_1, \cdots , X_k) \stackrel{def}{=} 
D_1 ^{(k)}[cr_2 F(- \vee \cdots \vee -, -)](X_1, \cdots ,X_{k-1}, X_k).
\end{equation}
To compute $A_k$, we expand the above second cross effect.
\begin{eqnarray}
cr_2 F(X_1 \vee \cdots \vee X_{k-1}, X_k) 
& \cong &
F(X_1 \vee \cdots \vee X_{k-1}) \otimes F(X_k) \nonumber\\
&\oplus& 
F( \Sigma^{-n+1} U(X_1 \vee \cdots \vee X_{k-1}) \otimes U(X_k)) \nonumber \\
&\oplus& 
F( \Sigma^{-n+1} U(X_1 \vee \cdots \vee X_{k-1}) \otimes U(X_k)) \nonumber \\
&&\otimes
[F(X_1 \vee \cdots \vee X_{k-1}) \oplus F(X_k) \nonumber \\ 
&&\oplus
F(X_1 \vee \cdots \vee X_{k-1}) \otimes F(X_k)] \nonumber
\end{eqnarray}
We compute $D_1 ^{(k)}$ of the three summands above. 

\begin{lemma}
\label{lemma:f1}
For $k \geq 3$, 
\begin{equation*}
D_1 ^{(k)} F(X_1 \vee \cdots \vee X_{k-1}) \otimes F(X_k) \cong 
A_{k-1}(X_1, \cdots ,X_{k-1}) \otimes D_1 F(X_k),
\end{equation*}
with $A_2 \cong D_1 F(X_1) \otimes D_1 F(X_2) \oplus 
\Sigma ^{n-1} D_1 F(X_1) \otimes D_1 F(X_2) \otimes D_1 F(\K)$, 
as computed above. 

\end{lemma}

\begin{proof}
The proof of this lemma follows from the following straight forward sequence of identities. 

\begin{eqnarray}
D_1 ^{(k)}F(X_1\vee \cdots \vee X_{k-1}) \otimes F(X_k)
& \cong & 
D_1 ^{(k-1)}F(- \vee \cdots \vee -)(X_1, \cdots ,X_{k-1}) \otimes D_1 F(X_k) \nonumber \\
& \cong &
D_1 ^{(k-1)} cr_2 F(X_1 \vee \cdots \vee X_{k-2}, X_{k-1}) \otimes D_1 F(X_k) 
\nonumber \\
& \cong &
A_{k-1} (X_1, \cdots ,X_{k-1}) \otimes D_1 F(X_k). \nonumber
\end{eqnarray}

\end{proof}

\begin{lemma}
\label{lemma:f2}
For all $k>2$,
 
\begin{eqnarray}
D_1 ^{(k)} F( \Sigma^{-n+1} U(-\vee \cdots \vee -) \otimes U(-))(X_1, \cdots , X_k) \cong
\nonumber \\
\bigoplus _{\sigma \in \Sigma_{k-1}} 
D_1 F(X_{\sigma(1)}) \otimes \cdots \otimes D_1 F(X_{\sigma(k-1)}) \otimes D_1 F(X_k)
\otimes \Sigma^{(k-1)(n-1)} D_1 F(\K).
\nonumber
\end{eqnarray}
\end{lemma}

\begin{proof}
Since it does not matter in which order we apply $D_1$'s, we begin by computing $D_1 ^k$,
i.e. the derivative with respect to the variable $X_k$. 
\begin{eqnarray}
\lefteqn{
D_1 ^{k} F( \Sigma^{-n+1} U(X_1\vee \cdots \vee X_{k-1}) \otimes U(-))(X_k)} \nonumber \\
& & \cong 
D_1 F( \Sigma^{-n+1} U(X_1\vee \cdots \vee X_{k-1}) \otimes D_1 U(-)(X_k)) \nonumber \\
& & \cong 
\Sigma^{-n+1} U(X_1\vee \cdots \vee X_{k-1}) \otimes \Sigma^{n-1} D_1 F(X_k) \otimes D_1 F(\K)
\nonumber \\
& & \cong 
U(X_1\vee \cdots \vee X_{k-1}) \otimes D_1 F(X_k) \otimes D_1 F(\K). \nonumber
\end{eqnarray}
Hence, recalling that $U$ is a functor of type $\frac{x}{1-x}$ and using our computations for 
these functors, we get 
\begin{eqnarray}
\lefteqn{
D_1 ^{(k)} F( \Sigma^{-n+1} U(X_1\vee \cdots \vee X_{k-1}) \otimes U(X_k))} \nonumber \\
& & \cong  
D_1 ^{(k-1)} U(X_1\vee \cdots \vee X_{k-1}) \otimes D_1 F(X_k) \otimes D_1 F(\K)
\nonumber \\
& & \cong  
\bigoplus _{\sigma \in \Sigma_{k-1}} 
D_1 U(X_{\sigma(1)}) \otimes \cdots \otimes U(X_{\sigma(k-1)}) \otimes D_1 F(X_k) \otimes D_1 F(\K)
\nonumber \\
& & \cong  
\bigoplus _{\sigma \in \Sigma_{k-1}} 
D_1 F(X_{\sigma(1)}) \otimes \cdots \otimes F(X_{\sigma(k-1)}) \otimes D_1 F(X_k) 
\otimes \Sigma^{(k-1)(n-1)} D_1 F(\K). 
\nonumber
\end{eqnarray}
\end{proof}

The final piece in which we are interested is 
\begin{eqnarray}
\lefteqn{
D_1 ^{(k)} F( \Sigma^{-n+1} U(X_1\vee \cdots \vee X_{k-1}) \otimes U(X_k))} \nonumber \\
& \otimes [F(X_1\vee \cdots \vee X_{k-1}) \oplus F(X_k) \oplus 
F(X_1\vee \cdots \vee X_{k-1}) \otimes F(X_k).
\nonumber
\end{eqnarray}
In other words we need to compute $D_1 ^{(k)}$ of the three summands above. Observe that 
the second and and the third summands are at least $2$-reduced in the variable $X_k$. 
Hence, $D_1 ^k$ applied to these summands produces objects equivalent to $0$. Consequently,
it remains to compute only 
$$D_1 ^{(k)} [F( \Sigma^{-n+1} U(X_1\vee \cdots \vee X_{k-1}) \otimes U(X_k)) 
\otimes F(X_1\vee \cdots \vee X_{k-1})].$$
We begin by applying $D_1 ^k$ to get 
\begin{eqnarray}
\Sigma^{-n+1} U(X_1\vee \cdots \vee X_{k-1}) \otimes D_1 U(X_k) \otimes D_1 F(\K) 
\otimes F(X_1\vee \cdots \vee X_{k-1}) \nonumber \\
\cong 
U(X_1\vee \cdots \vee X_{k-1}) \otimes F(X_1\vee \cdots \vee X_{k-1}) 
\otimes D_1 F(X_k) \otimes D_1 F(\K)
\nonumber
\end{eqnarray}
Thus, we need to understand
$$D_1 ^{(k-1)}[U(X_1\vee \cdots \vee X_{k-1}) \otimes F(X_1\vee \cdots \vee X_{k-1})].$$
Recall that $U$ is a functor of type $\frac{x}{1-x}$ and hence 
$U(X_1\vee \cdots \vee X_{k-1})$ can be expanded as a sum. Note that the only terms 
of $U(X_1\vee \cdots \vee X_{k-1})$ which are not at least $2$-reduced in at least 
one of the variables are
$$ U(X_1) \oplus \cdots \oplus U(X_{k-1}) \oplus
\bigoplus_{i \neq j} U(X_i) \otimes U(X_j) \oplus \cdots \oplus
\bigoplus_{\sigma \in \Sigma_{k-1}} U(X_{\sigma(1)}) \otimes U(X_{\sigma(k-1)}).$$
A typical summand in the above sum is of the form 
$$U(X_{i_1}) \otimes \cdots \otimes U(X_{i_s}),$$
with $s \leq k-1$ and ${i_l} \neq {i_t}$ for $l \neq t$. 
Consequently, it is enough to understand 
$$D_1 ^{(k-1)} [U(X_{i_1}) \otimes \cdots \otimes U(X_{i_s}) \otimes 
F(X_1\vee \cdots \vee X_{k-1})].$$
To simplify the indexing, let us assume that in the above term, $i_j=j$, in other 
words, we compute
\begin{eqnarray}
\lefteqn{
D_1 ^{(k-1)} [U(X_{1}) \otimes \cdots \otimes U(X_{s}) \otimes 
F(X_1 \vee \cdots \vee X_{k-1})]} \nonumber \\
&& \cong D_1 ^{(k-1)} [U(X_{1}) \otimes \cdots \otimes U(X_{s}) \otimes 
(F(X_1) \oplus F(X_2 \vee \cdots \vee X_{k-1}) \nonumber \\
&& \oplus F( \Sigma^{-n+1} U(X_1) \otimes F(X_2 \vee \cdots \vee X_{k-1}) \nonumber \\
&& \oplus F( \Sigma^{-n+1} U(X_1) \otimes F(X_2 \vee \cdots \vee X_{k-1})) \otimes (\cdots))]
\nonumber \\
&& \cong D_1 ^{(k-1)} [U(X_{1}) \otimes \cdots \otimes U(X_{s}) 
\otimes F(X_2 \vee \cdots \vee X_{k-1})] \nonumber
\end{eqnarray}
since all the other terms are at least $2$-reduced in $X_1$. Iterating this process for 
$X_2, \cdots , X_s$, we get that 
\begin{eqnarray*}
\lefteqn{
D_1 ^{(k-1)} [U(X_{1}) \otimes \cdots \otimes U(X_{s}) \otimes 
F(X_1 \vee \cdots \vee X_{k-1})]} \\   
&&\cong 
D_1 ^{(k-1)} [U(X_{1}) \otimes \cdots \otimes U(X_{s}) \otimes 
F(X_{s+1} \vee \cdots \vee X_{k-1})] \\   
&&\cong
D_1 U(X_1) \otimes \cdots \otimes D_1 U(X_{s}) \otimes A_{k-s-1}(X_{s+1}, \cdots , X_{k-1}) \\  
&&\cong
\Sigma^{s(n-1)} 
D_1 F(X_1) \otimes \cdots \otimes D_1 F(X_{s}) \otimes A_{k-s-1}(X_{s+1}, \cdots ,X_{k-1}),
\end{eqnarray*}
where $A_{k-s-1}$ is as it is defined in Equation~\ref{eq:A_k}.

We combine these computations in the following lemma by first agreeing to augment 
the Equation~\ref{eq:A_k} by defining $A_1(X) \stackrel{def}{=} D_1 F(X)$.
\begin{lemma}
\label{lemma:f3}
For all $k>2$ and $X_1 = \cdots = X_k=X$, and abbreviating $A_{k-s-1}(X, \cdots ,X)$ to 
$A_{k-s-1}$, we have  
\begin{eqnarray*}
\lefteqn{
D_1 ^{(k)} F( \Sigma^{-n+1} U(X \vee \cdots \vee X ) \otimes U(X ))}  \\
&& \otimes [F(X \vee \cdots \vee X ) \oplus F(X) \oplus 
F(X \vee \cdots \vee X) \otimes F(X)] \\
&& \cong \bigoplus_{s=1} ^{k-2} \Sigma^{s(n-1)} 
[\bigoplus_{\sigma \in P(k-1,s)} \bigotimes_1 ^s D_1 F(X) 
\otimes A_{k-s-1}] 
\otimes D_1 F(X) \otimes D_1 F(\K),
\end{eqnarray*}
where $P(k-1,s)$ is the set of $s$ different ordered letters from a collection 
of $k-1$ letters.
\end{lemma}
The above Lemmas~\ref{lemma:f1},~\ref{lemma:f2},~\ref{lemma:f3} give a complete description
of the differentials $\bigtriangledown ^k F(X, \cdots , X; 0)$, which are of course the 
essential ingredients for computing the layers by Corollary~\ref{Randy2}. We agree to 
augment the Equation~\ref{eq:A_k} further by setting $A_0 = \K$ and combine these
lemmas to obtain the following proposition.
\begin{proposition}
\label{prop:diff}
The differentials of functors of type $f_n$ are as follows.

\noindent 
For k=2,
\begin{eqnarray*}
\lefteqn{
\bigtriangledown ^2 F(X, X;0) \cong A_2 \cong}\\
&&D_1 F(X) \otimes D_1 F(X) \oplus 
\Sigma ^{n-1} D_1 F(X) \otimes D_1 F(X) \otimes D_1 F(\K).
\end{eqnarray*}
For $k>2$,
\begin{eqnarray*}
\lefteqn{
\bigtriangledown ^k F(X, \cdots , X; 0) \cong A_k \cong} \\
&& A_{k-1} \otimes D_1 F(X) \oplus \\
&& \bigoplus_{s=1} ^{k-1} \Sigma^{s(n-1)} 
[\bigoplus_{\sigma \in P(k-1,s)} [\bigotimes_1 ^s D_1 F(X)] 
\otimes A_{k-s-1}] 
\otimes D_1 F(X) \otimes D_1 F(\K), 
\end{eqnarray*}
where, as before, $P(k-1,s)$ is the set of $s$ different ordered letters from a 
collection of $k-1$ letters, while $A_0 =\K$ and $A_1= D_1 F(X)$.
\end{proposition}

\section{Functional Equations and Operads}

We have already commented on connection between functors of type 
$f_n$ and the homology of little $n$-cubes operads, which was suggested by the fact that
they have the same set of associated algebras. 
In fact, it is in part with this relationship in mind that we considered these functors.
The following discussion is intended to explore this relationship further.

As an easy consequence of  the defining equation of functors of type $f_n$ 
(Definition~\ref{def:f_n}), we note that (with an additional hypothesis that 
$F$ is continuous) as $n$ increases, so does the connectivity of the term
$$F(\Sigma^{-n+1} U X \otimes U Y) \cong 
F(\Sigma^{-n+1} \P_1\Sigma^{n-1} F(X) \otimes \P_1\Sigma^{n-1} F(Y)).$$
Thus, the equation degenerates
to the defining equation of exponential functors. This, of course, is in agreement with 
the fact that as a function of real variables $e^x -1 $ is the limit of the sequence of
functions $\{f_n\}$. 
Moreover, by analyzing the differentials $D_n$ of functors of type $f_n$ that we 
just computed, it is evident that the connectivity of the terms 
$$\bigoplus_{s=1} ^{k-1} \Sigma^{s(n-1)} 
[\bigoplus_{\sigma \in P(k-1,s)} [\bigotimes_1 ^s D_1 F(X)] 
\otimes A_{k-s-1}] 
\otimes D_1 F(X) \otimes D_1 F(\K)$$
of $D_n$ increases with $n$, and hence in the limit case, the layers are the same as the 
layers of exponential functors.

However, if we set $n=1$ in the equation of Definition~\ref{def:f_n}, 
we get functors that factor through the category of usual Poisson algebras, while the 
functors of type $\frac{x}{1-x}$ (Definition~\ref{def:ass}) factor through the category 
of associative algebras, and of course, the associative and Poisson algebras are quite 
different. A similar situation occurs in ~\cite{Getz} where E. Getzler and J.D.S. Jones 
prove that the homology little $n$-cubes operads ${\bf e_n}$ are isomorphic to the Poisson
operads ${\bf p_n}$ for $n>1$ (Theorem 1.6, ~\cite{Getz}), while ${\bf e_1}$ is not 
isomorphic to ${\bf p_1}$. However, ${\bf e_1}$ is naturally a filtered operad, and the 
associated graded operad
is the Poisson operad  ${\bf p_1}$ (Section 5.3, ~\cite{Getz}). We have a similar result. 
To make this statement precise, let us consider the functors into $Ch(\K)$ of type 
$\frac{x}{1-x}$ again. We already showed that if $F$ is such a functors then $F(X)$ is
naturally an associative algebra. In fact, it is equipped with an even richer algebraic 
structure if we restrict ourselves to functors whose domain is $Ch(\K)$. Using the fact 
that finite coproducts are also products in $Ch(\K)$ we define a comultiplication map  
 $\Delta_{F(X)}: F(X) \rightarrow F(X) \otimes F(X)$ for such functors:

\hspace{25pt} $ F(X) \stackrel{\Delta}{\rightarrow} F(X \vee X) \rightarrow F(X) \otimes F(X) 
  \oplus F(X) \otimes F(X) \stackrel{+}{\rightarrow} F(X) \otimes F(X) .$

\noindent
Combined with the associative multiplication defined earlier,
$\Delta_{F(X)}$ makes $F(X)$ into a Hopf algebra. In addition, from the symmetry of the 
definition of $\Delta_{F(X)}$ it is immediate that $F(X)$ is cocomutative. Thus, $F$ factors
through the category $\CH$ of cocomutative Hopf algebras.     
Hence, by the structure theorem of Milnor and Moor ~\cite{Milnor}, $F(X)$ is isomorphic
to the universal enveloping algebra $U(L)$ of some Lie algebra $L$. In fact, this can be
done functorially to produce a commutative diagram

$$
\xymatrix{
Ch(\K) \ar[r]^-{\hat{F}} \ar[d]^-{\tilde{F}}
&
\A \\
\CH \ar[r]^-{P} 
&
\L \ar[u]_-{U}
}
$$
where $\hat{F}$ exists by Lemma~\ref{lemma:ass}, $\tilde{F}$ by the argument above, and $P$ is the 
inverse functor of the universal enveloping algebra functor $U$ considered as a functor 
into $\CH$ (Theorem 5.18 of ~\cite{Milnor}). Here, as before, $\L$ is the category of Lie algebras
over the field $\K$ of characteristic $0$. 

Now we define a filtration on $F(X)$ which 
is a essentially the {\textit {Lie filtration}} of~\cite{Milnor} (see Definition 5.12 
of~\cite{Milnor}). In other words, roughly speaking, it is the filtration of $F(X)$ given 
by powers of $P \circ \tilde{F}(X)$,
thought of as a subobject of $\hat{F}(X)$ via the universal enveloping algebra 
functor, i.e. 
$F_p(X) = \sum_{q \leq p} [{P \circ \tilde{F}(X)}]^q$. 
More precisely,

$$F_1(X) = Im[\K \oplus P \circ \tilde{F}(X) \rightarrow \hat{F}(X)]$$ 
$$F_{p+1}(X) = Im [F_p(X) \oplus ((P \circ \tilde{F}(X)) \otimes F_p(X)) \rightarrow \hat{F}(X)].$$

\noindent
The associated graded object is 
defined to be $A(X)=F_1 \oplus F_2/F_1 \oplus \cdots F_p/F_{p-1} \oplus \cdots$. Now we 
state the promised result.

\begin{proposition}

If $F: Ch(\K) \rightarrow Ch(\K)$ is a functor of type $\frac{x}{1-x}$ then the 
correspondence $X \longmapsto A(X)$ defines a functor of type $f_1$.

\end{proposition}

\begin{proof}
Throughout this proof we omit the forgetful functors where no confusion may arise 
to ease the notation. Thus, for example, when we refer to a functor $F$ into Lie algebras 
as a logarithmic functor, naturally we mean $F$ composed with the forgetful functor 
$U_{log}$ from Lie algebras to chain complexes, as by definition, all our functors 
of different types have the category $Ch(\K)$ as their target.
   
First observe that $A(X)$ is a functor, i.e. that 
any map $f:X \rightarrow Y$ in $Ch(\K)$ induces a morphism of filtered 
algebras. Indeed, as $\tilde{F}$ and $P$ are functors we have a map 
$\T \circ P \circ \tilde{F}(X) \rightarrow \T \circ P \circ \tilde{F}(Y)$, which takes 
components of degree $\leq q$ to components of degree $\leq q$. 
Recall that the universal enveloping functor $U$ on a Lie algebra ${\g}$
can be constructed as a quotient of $\T \g$ by the ideal generated by relations 
$i([x,y])=i(x)i(y) - i(y)i(x)$, where $i:\g \to \T \g$ is the obvious embedding,
and $x,y \in \g$.   
Hence we have that $f^*(F_p(X)) \subseteq F_p(Y)$, 
where $f^*$ is the map induced by $f$. Thus, $f$ produces a morphism of filtered algebras 
which in its turn induces a map of corresponding graded objects $A(X) \rightarrow A(Y)$ 
making $A$ into a functor. 

It is a well known fact about Lie algebras that the graded algebra associated to the 
{\textit {Lie filtration}} (that we just presented) is graded commutative 
(see Section 7.3 of~\cite{Weib} for example). 
Moreover, if we denote by $S$ the functor that takes the module underlying 
the Lie algebra $\g$ to the chain 
complex underlying the graded algebra associated with the {\textit {Lie filtration}} of 
$U(\g)$, then it satisfies the equation 
$S(\g \oplus \h) \cong S(\g) \oplus S(\h) \oplus S(\g) \otimes S(\h)$. 
where $\g \oplus \h$ is the coproduct of $\g$ and $\h$ in $Ch(\K)$, with the Lie bracket 
defined in the obvious way: $[(a,b), (\bar{a},\bar{b})]=([a,\bar{a}]_{\g},[b,\bar{b}]_{\h})$.
This is the case since in characteristic $0$, the graded commutative algebra
associated to the {\textit {Lie filtration}} of $U(\g)$ 
is isomorphic to the free symmetric algebra on $\g$.
Thus, $S$ is an exponential functor.

To complete the proof of the proposition we simply observe that the functor $A$ is the 
composition of $S$ with $P \circ \tilde{F}$, and since  $P \circ \tilde{F}$ is a coproduct 
preserving functor into Lie algebras ($P$ is coproduct preserving by Theorem 5.9 of 
~\cite{Milnor}) and thus logarithmic, the computations done in                  
Example~\ref{example:f_n} for $n=1$ imply that $A=S \circ P \circ \tilde{F}$ satisfies 
the defining equation of functors of type $f_1$.

\end{proof}

We conclude by considering the example of the forgetful functor from the category 
of $n$-Poisson algebras in greater detail. Its relevance to the algebras over the homology 
of little $n$-cubes operads, as well as to the homology of configuration spaces,
is the main reason for our interest. 

Recall that the triple associated with the category of $n$-Poisson algebras is of the 
form
\begin{equation}
\label{eq:trip}
T_n (X) = \bigoplus _k {\bf p_n}(k) \otimes_{\Sigma_k} X^{\otimes k},
\end{equation} 
where ${\bf p_n}$ is the $n$-Poisson operad. The algebras over this triple are precisely 
the $n$-Poisson algebras. Of course, for $n \geq 2$ the operads ${\bf p_n}$ are 
equivalent to the homology of little $n$-cubes operads ${\bf e_n}$ (e.g. see ~\cite{Getz}). 
Note that we can read the operad off the triple (~\ref{eq:trip}), 
by computing the differentials evaluated at the ground field $\K$ viewed as a chain 
complex concentrated in degree $0$:
$$\bigtriangledown ^k F(\K, \cdots , \K; 0).$$
Observing that $D_1 T_n (\K) \cong \K$, we can rewrite the Proposition~\ref{prop:diff}
as follows:

\noindent
For $k=2$,
\begin{equation}
\label{eq:t_n1}
\bigtriangledown ^2 T_n (\K,\K;0) \cong A_2 \cong \K \oplus \Sigma ^{n-1} \K.  
\end{equation}

\noindent
For $k>2$,
\begin{equation}
\label{eq:t_n2}
\bigtriangledown ^k T_n (\K, \cdots , \K;0) \cong A_k \cong A_{k-1} \oplus
\bigoplus _{s=1} ^{k-1} \Sigma^{s(n-1)}
[\bigoplus_{\sigma \in P(k-1,s)} A_{k-s-1}],
\end{equation}
where $A_i$ and $P(k-1,s)$ are as before. In fact, we can provide a better description.

\begin{lemma}
\label{lem:difhom}
For $k \geq 2$,
\begin{equation}
\label{eq:homsph}
\bigtriangledown ^k T_n (\K, \cdots , \K;0) 
\cong \bigotimes _{j=1} ^{k-1} H^{\ast} (\bigvee _j S^{n-1}),
\end{equation}
where on the right hand side we have the tensor product of 
cohomologies of the wedge of $j$ copies of the sphere $S^{n-1}$.
Here and in what follows the homologies are taken with coefficients in the field $\K$. 
\end{lemma}

\begin{proof}
We prove by induction on $k$. Since we intend to use the formulas (~\ref{eq:t_n1}) and 
(~\ref{eq:t_n2}), we have to consider both $k=2$ and $k=3$ for the base case. 

For $k=2$, we need to show that 
$\bigtriangledown ^2 T_n (\K, \K;0) \cong H^{\ast}(S^{n-1})$. This is immediate from 
Equation~\ref{eq:t_n1}.

For $k=3$, note that 
$\bigotimes _{j=1} ^{2} H^{\ast} (\bigvee _j S^{n-1}) \cong H^{\ast}(S^{n-1}) \otimes 
H^{\ast}(S^{n-1} \vee S^{n-1})$. Hence by Mayer-Vietoris (e.g. see~\cite{Bredon}), the 
right hand side of Equation~\ref{eq:homsph} has a copy of $\K$ in degree $0$, three 
copies of $\K$ in degree $n-1$ and two copies of $\K$ in degree $2(n-1)$. 
On the other hand by Equations~\ref{eq:t_n1} and~\ref{eq:t_n2},
\begin{equation*}
\bigtriangledown ^3 T_n (\K,\K,\K;0) \cong \K \oplus \Sigma ^{n-1} \K \oplus
\Sigma ^{n-1}(\K \oplus \K) \oplus \Sigma ^{2(n-1)}(\K \oplus \K), 
\end{equation*}
which proves the case $k=3$.

Suppose lemma holds for $k=l-1$. For $k=l$, by Mayer-Vietoris, on the right hand side of 
Equation~\ref{eq:homsph} we have
\begin{eqnarray*}
\lefteqn{
\bigotimes _{j=1} ^{l-1} H^{\ast} (\bigvee _j S^{n-1})} \\
&& \cong 
\bigotimes _{j=1} ^{l-2} H^{\ast} (\bigvee _j S^{n-1}) \otimes 
H^{\ast} (\bigvee _{l-1} S^{n-1}) 
\cong 
\bigotimes _{j=1} ^{l-2} H^{\ast} (\bigvee _j S^{n-1}) \otimes 
[\K \oplus \Sigma ^{n-1} \K ^{\oplus {l-1}}]\\
&& \cong A_{l-1} \otimes[\K \oplus \Sigma ^{n-1} \K^{\otimes l-1}] 
\cong A_{l-1} \otimes \K \oplus _{l-1} \otimes \Sigma ^{n-1} \K^{\otimes l-1} \\
&& \cong A_{l-1} \oplus [A_{l-2} \oplus \bigoplus _{s=1} ^{l-2} \Sigma^{s(n-1)}
(\bigoplus_{P(l-2,s)} A_{l-s-2})] \otimes \Sigma ^{n-1} \K^{\otimes l-1}\\
&& \cong A_{l-1} \oplus A_{l-2} \otimes \Sigma ^{n-1} \K^{\otimes l-1} \oplus 
\bigoplus _{s=1} ^{l-2} \Sigma^{(s+1)(n-1)} [\bigoplus_{P(l-2,s)} A_{l-s-2} 
\otimes \K^{\otimes l-1}].
\end{eqnarray*}
We renumber the sum in the last equation by setting $s+1=t$ to get 
\begin{eqnarray*}
\lefteqn{
\bigotimes _{j=1} ^{l-1} H^{\ast} (\bigvee _j S^{n-1})} \\
&& \cong A_{l-1} \oplus A_{l-2} \otimes \Sigma ^{n-1} \K^{\otimes l-1} \oplus
\bigoplus _{t=2} ^{l-1} \Sigma^{t(n-1)} [\bigoplus_{P(l-2,t-1)} A_{l-t-1} 
\otimes \K^{\otimes l-1}].
\end{eqnarray*}
We abuse the notation to denote by $P(l-2,t-1)$ the number of elements in $P(l-2,t-1)$. 
Then 
$$P(l-2,t-1) (l-1) = \frac{(l-2)!}{(l-t-1)!} (l-1) = \frac{(l-1)!}{(l-t-1)!} = P(l-1,t).$$
These identities allow us to simplify our equations further:
\begin{eqnarray*}
\lefteqn{
\bigotimes _{j=1} ^{l-1} H^{\ast} (\bigvee _j S^{n-1})} \\
&& \cong A_{l-1} \oplus A_{l-2} \otimes \Sigma ^{n-1} \K^{\otimes l-1} \oplus
\bigoplus _{t=2} ^{l-1} \Sigma^{t(n-1)} [\bigoplus_{P(l-1,t)} A_{l-t-1}]\\
&& \cong  A_{l-1} \oplus
\bigoplus _{t=1} ^{l-1} \Sigma^{t(n-1)} [\bigoplus_{P(l-1,t)} A_{l-t-1}]
\cong A_l \cong \bigtriangledown ^l T_n (\K, \cdots , \K;0),
\end{eqnarray*}
thus completing the proof of the lemma.

\end{proof}

Recall that one of the descriptions of the homology of little $n$-cubes operads is via 
configuration spaces (see Section~\ref{sect:operad}), which is the reason why our 
Lemma~\ref{lem:difhom} is reminiscent of results of F.Cohen (see ~\cite{Cohen}), 
who applied a theorem by 
Fadell and Neuwirth (see~\cite{FadNeu}) to describe the homology of configuration 
spaces of embedding into $\R ^n$:

\begin{lemma}[Cohen]
Additively,
\begin{equation*}
H^ \ast (F_n(k)) \cong  \bigotimes _{j=1} ^{k-1} H^{\ast} (\bigvee _j S^{n-1}),
\end{equation*}    
where $F_n(k)$ is the configuration space of embedding into $\R ^n$ as it is defined 
in Section~\ref{sect:operad}.
\end{lemma}

Consequently, these computations of homology of configuration spaces can be recovered from 
Lemma~\ref{lem:difhom}. Similarly, our techniques can be used to reproduce the description 
of the action of symmetric groups on $H^ \ast (F_n(k))$ (Section 7 of~\cite{Cohen}).

\end{document}